\documentstyle[12pt,leqno,amstex,amssymb,amsfonts,eucal]{article}
\begin{document}
\oddsidemargin= 12mm
\topmargin= 35mm
\textwidth=170mm
\textheight=270mm
\pagestyle{plain}
\newcounter{r}
\newcounter{eqh}
\newcommand{\Ker}{Ker}
\newcommand{\sign}{sign}
\newcommand{\Mat}{Mat}
\renewcommand {\Re}{Re}
\renewcommand{\Im}{Im}
\begin{center}
\large
{\bf  ON THE DIFFERENCE EQUATION OF THE POINCAR\'E TYPE }\\
\vskip 10pt
{\bf L.A. Gutnik}
\end{center}
\vskip5pt
{\hfill \sl Dedicated to the memory of Professor A.O. Gelfond.}
\vskip10pt
\small Abstract. In the paper are proved theorems, which
amplify the results of my paper
"On the difference equation
of Poincar\'e type (Part 3)",
Max-Plank-Institut f\"ur Mathematik, Bonn,
 Preprint Series, 2004, 09, 1 -- 34.
\vskip 10pt
\normalsize
\vskip.10pt
\begin{center}{\large\bf Table of contents}\end{center}
\noindent \S 0. Foreword.
\vskip.4pt
\noindent \S 1. Some preparatory results.
\vskip.4pt
\noindent \S 2. Proof of the theorem 8.
\vskip.4pt
\noindent \S 3. Proof of the Theorem 9.
\vskip.4pt
\noindent \S 4. Proof of the Theorem 10.
\vskip.4pt
\noindent \S 5. The case of the general
differrence equation of the Poincar\'e type.
{\begin{center}\large\bf\S 0. Foreword.\end{center}}
In [\ref{r:ci0}]: was proved the following

\bf Theorem 7. \it Let $s\in{\Bbb N}-1,n\in{\Bbb N,}$
$$
a^\sim_i\in{\Bbb C}, \; a_i(\nu)\in{\Bbb C},
$$
\begin{equation}
\label{eq:zc}
 a_n(\nu)=1,\; a_i(\nu)-a^\sim_i=O(1/(\nu+1))
\end{equation}
for $\nu\in{\Bbb N}-1$ and $i=0,\,\ldots,\,n.$
Let us consider the following difference equation
\begin{equation}
\label{eq:za0}
\sum\limits_{k=0}^na_k(\nu)y(\nu+k)=0,
\end{equation}
where $\nu\in {\Bbb N}-1.$
\newpage
\pagestyle{headings}
\topmargin= -15mm
\textheight=250mm
\markright{\footnotesize\bf L.A.Gutnik,
  ON THE DIFFERENCE EQUATION OF THE POINCAR\'E TYPE}
 For $m\in{\Bbb N}$ let $V_m$ denotes the linear
 over ${\Bbb C}$ space of solutions  $y=y(\nu)$ of the equation
\begin{equation}\label{eq:zb}
\sum\limits_{k=0}^na_k(\nu)y(\nu+k)=0,
\end{equation}
where $\nu\in m+{\Bbb N}-1.$
Let the absolute values of all the roots of
 the characteristical polynomial
\begin{equation}\label{eq:zd}
T(z)=\sum\limits_{k=0}^na_k^\sim z^k
\end{equation}
\noindent are among the numbers
$\{\rho_i\colon1\leq i\leq 1+s\}$ such that
$\rho_{s+1}=0$ and $\rho_j<\rho_i$ for $1\leq i< j\leq s+1.$
Let $e_i$ and $k_i$ denote respectively the sum and the maximum
 of the multiplicities of those roots,
 whose absolute value is equal to the number $\rho_i,$
 where $i=1,\,\ldots,\,s+1,$
 and let $k^\ast=k_{s+1}.$ We suppose that, if $s>0,$ then
\begin{equation}\label{eq:zd1}
e_i>0
\end{equation}
for $i=1\,\ldots,\,s.$
 For given $y=y(\nu)$ in ${\Bbb C}^{m-1+{\Bbb N}}$, let
$$\omega_{n,y}(\nu) =
\max (\vert y(\nu)\vert\,,\ldots\,,\vert y(\nu+n-1)\vert).$$
Then there exist $A>0,\,m\in{\Bbb N},\,\alpha^\wedge(\nu)>0$
 with $\nu \in{m+{\Bbb N}-1}$ and
 the subspaces $V_{m,1}^\vee,\,\ldots,\,V_{m,s+1}^\vee$ such that
$$\lim\limits_{\nu\to\infty}\alpha^\wedge(\nu)=0,$$
$$V_m=V_{m,1}^\vee\oplus\,\ldots,\oplus\,V_{m,s+1}^\vee,\,
\dim_{\Bbb C}(V_{m,i}^\vee)=e_i,\;1\leq i\leq s+1,$$
and, if $y\in V_{m,\theta}^\vee$ for some $\theta\in\{1,...,s\},$ then
\begin{equation}\label{eq:za1}
\exp(-A(ln(\nu)+\nu^{1-1/k_\theta}))(\rho_\theta)^{\nu}\omega_n(y)(m)
\le
\omega_{n,y}(\nu)\end{equation}
for $\nu\in m+{\Bbb N}-1;$
 moreover, the spaces
$$V_{m,j}^\wedge = V_{m,j}^\vee \oplus\,\ldots\,\oplus V_{m,s+1}^\vee,$$
where $j=1\,\ldots,\, s+1,$ and, if $s\ge1$  natural projections $\pi_j,$
$$V_{m,j}^\wedge\mapsto V_{m,j}^\vee,$$
 where $j=1\,\ldots,\,s,$ have the folloving properties:
\vskip2pt
if $y\in V_{m,\theta}^\wedge$ for some $\theta\in\{1,...,s\}$, then
\begin{equation}\label{eq:za2}
\omega_{n,y}(\nu)\le
\exp(A(\ln(\nu)+\nu^{1-1/k_\theta}))(\rho_\theta)^{\nu}\omega_{n,y}(m),
\end{equation}
\newpage
\begin{equation}\label{eq:za4}
(\omega_{n\pi_{\theta(y)}}(m)-\alpha(\nu)\omega_{n,y}(m))
(\rho_\theta)^{\nu}
\end{equation}
$$\exp(-A(\ln(\nu)+\nu^{1-1/k_\theta}))
\le\omega_{n,y}(\nu),$$
where $\nu\in m+{\Bbb N}-1;$ if
\begin{equation}\label{eq:za3}
k^\ast>0,
\end{equation}
 and $y\in V_{m,s+1}^\vee\,(=V_{m,s+1}^\wedge)$, then
\begin{equation}\label{eq:xza3}
\vert y(\nu)\vert\le(A/\nu)^{\nu/k^\ast}\omega_{n,y}(m),
\end{equation}
where $\nu\in m+{\Bbb N}-1.$

\bf Remark 1. It follows from the Theorem 7 that the space
 $V_{m,\theta}^\wedge,$ where $\theta=1,\,\ldots,\,s+1,$
 does not depend from the construction
 and is defined uniquelly by means of the equality
$$V_{m,\theta}^\wedge=\{y\in V_m\colon \limsup\vert y(\nu)\vert^{1/\nu}
\le\rho_\theta\}.$$ \rm

The presence of unkown $\alpha(\nu)$ (even tending to zero) in (\ref{eq:za4})
 constrict the possibilities of the application of this Theorem.
 Of course, in view of (\ref{eq:za1}), in the case $e_{s+1}=0,\,\theta=s$
 this $\alpha(\nu)$ cannot play devil with us ,
 because it vanishes then,
 but such happy case (see, for example, [\ref{r:ch0}]-[\ref{r:zch0}])
 is rather an exeption from the rule.
 One may attempt to estimate the specified value but
 it would be better to get rid from it at all. With this goal
 I prove here the following

\bf Theorem 10. \it Let are fulfilled all the conditions of the Theorem 7.
Let further $V$ be an arbitrary linear subspace of $V_m$ such that
$$V\cap V_{m,\theta+1}=\{0\},$$
where $\theta\in\{1\,\ldots,\, s\}.$

Then for this $V$ there exists a constant $A^\ast=A^\ast(V)>0$ such that
\begin{equation}\label{eq:za5}
\exp(-A^\ast(ln(\nu)+\nu^{1-1/k}))(\rho_\theta)^{\nu}\omega_n(y)(m)
\le
\omega_{n,y}(\nu)\end{equation}
where $y\in V,\,k=\max(k_1,\,\ldots\,k_s)$ and $\nu\in m+{\Bbb N}-1.$

\rm First I prove the vollowing

\bf Theorem 8. \it Let are fulfilled all the conditions of the Theorem 7.
Let further
$$V_{m,j}^\ast = V_{m,1}^\vee \oplus\,\ldots\,\oplus V_{m,j}^\vee,$$
where $j=1\,\ldots,\, s+1,$ (and $V_{m,1}^\ast=V_{m,1}^\vee.$)
Then vor $V=V_{m,\theta}^\ast$
with $\theta\in\{1,...,s\},$
holds the assertion of the Theorem 10.

\rm Then I prove

\bf Theorem 9. \it Let for some $\theta\in\{1\,\ldots,\, s\}$ is
 given a linear map $\xi_\theta$ of the space $V_{m,\theta}^\ast$
 into $V_{m,\theta+1}^\wedge.$ Let $I^\ast_\theta$ is
 the identity map $V_{m,\theta}^\ast\to V_{m,\theta}^\ast$ Then for
$$V=(I_\theta^\ast+\xi_\theta)(V_{m,\theta}^\ast)$$
holds the assertion of the Theorem 10.

\rm In the section 4 I prove the Theorem 10.

And in the section 5 I discuss, what will take place, if instead
(\ref{eq:zc}) the following conditions hold:
\begin{equation}\label{eq:xzc}
\lim\limits_{\nu\to\infty} a_i(\nu) = a^\sim_i,
\end{equation}
where $i=0,\,\ldots,\,n,$
\begin{equation}\label{eq:yzc}
 a_n(\nu)=1,
\end{equation}
where $\nu\in{\Bbb N}-1.$
{\begin{center}\large\bf \S 1. Some preparatory results.
\end{center}}

\bf Lemma 1. \it Let $a\in{\Bbb N},\,b\in{\Bbb N}-1+a,\,C>0$ Then
\begin {equation}\label{eq:a}
0<\sum\limits_{\kappa=a}^b\ln(1+C/\kappa)
\le\end{equation}
$$\ln(1+C/a)+b\ln(1+C/b)-a\ln(1+C/a)+
C\ln((b+C)/(a+C)).$$

\bf Proof. \rm See the proof of the Lemma 1 in [\ref{r:ci0}]. $\blacksquare$

\bf Corollary. \it If $a\in{\Bbb N},\,b\in{\Bbb N}+a-1,\,b<2a,\,C>0,$
then
$$\sum\limits_{\kappa=a}^b\ln(1+C/\kappa)\le3C.$$

\bf Proof. \rm See the Proof of the Corollary
 of the Lemma 1 in [\ref{r:ci0}]. $\blacksquare$

\bf Lemma 2. \it ([\ref{r:bj}], Lemma 2, [\ref{r:ad}], Lemma 2,
[\ref{r:cf0}], Lemma 8.)

Let $A\in Mat_n({\Bbb C})$ an let $k$ is a maximal order of
 its Jordan blocks. Then there exists a constante $\gamma^\ast(A)>0$
 with the following properties:

for any $\varepsilon>0$ there exists a norm $p_{A,\varepsilon}$
 on ${\Bbb C}^n$ such that
\begin {equation}\label{eq:ec}
 p_{A,\varepsilon}\le\gamma^\ast(A)(\max(1,1/\varepsilon)^{k-1}h,
\end{equation}
\begin {equation} \label{eq:ed}
h\le\gamma^\ast(A)(\max(1,\varepsilon)^{k-1}p_{A,\varepsilon},
\end{equation}
\begin {equation} \label{eq:ee}
(p_{A,\varepsilon})^\sim\le(\gamma^\ast(A))^2
(\max(\varepsilon,1/\varepsilon)^{k-1}h^\sim,
\end{equation}
\begin {equation} \label{eq:ef}
h^\sim\le(\gamma^\ast(A))^2
(\max(\varepsilon,1/\varepsilon)^{k-1}(p_{A,\varepsilon})^\sim,
\end{equation}
\begin {equation} \label{eq:eg}
\Vert A\Vert_{sp}\le
(p_{A,\varepsilon})^\sim\le \Vert A\Vert_{sp}+(\sign(k-1))\varepsilon,
\end{equation}
where $\Vert A\Vert_{sp}$ denotes the maximum of the absolute values
 of eigenvalues of the matrix $A.$ If, moreover,
\begin {equation} \label{eq:eh}
\det(A)\ne0,\,
 \left\Vert A^{-1}\right\Vert_{sp}^{-1}>(\sign(k-1))\varepsilon,
\end{equation}
then
\begin {equation} \label{eq:ei}
 \left\Vert A^{-1}\right\Vert_{sp}\le
\end{equation}
$$(p_{A,\varepsilon})^\sim(A^{-1})\le
\left( \left\Vert A^{-1}\right\Vert_{sp}^{-1}-
(\sign(k-1))\varepsilon\right)^{-1}=$$
$$(\sign(k-1))\left\Vert A^{-1}\right\Vert_{sp}+
\varepsilon\left\Vert A^{-1}\right\Vert_{sp}
\left( \left\Vert A^{-1}\right\Vert_{sp}^{-1}-
(\sign(k-1))\varepsilon\right)^{-1}.$$
\bf Proof. \rm See the proof of the Lemma 8 in [\ref{r:cf0}]. $\blacksquare$

\bf Corollary. \it If all the eigenvalues of the matrix $A$ are simple, then
\begin {equation}\label{eq:fh}
(p_{A,\varepsilon})^\sim=\Vert A\Vert_{sp}.
\end{equation}
If, moreover,
\begin {equation} \label{eq:fi}
\det(A)\ne0,
\end{equation}
then
\begin {equation}
\label{eq:gj}
(p_{A,\varepsilon})^\sim(A^{-1})=
\left(\left\Vert A^{-1}\right\Vert_{sp}\right)^{-1}.
\end{equation}

\bf Proof. \rm See the proof of the Corollary
of the Lemma 8 in [\ref{r:cf0}]. $\blacksquare$

\bf Lemma 3. \it Let all the conditions of the Theorem 7 are fulfilled,
 and let  $$k=\max(k_1,\,\ldots,\,k_s).$$
 Then there exist $A>0,\,m\in{\Bbb N}$ such that
\begin{equation}\label{eq:xza2}
\omega_{n,y}(\nu)\le\exp(A(ln(\nu)+\nu^{1-1/k}))(\rho_1)^{\nu}\omega_n(y)(m)
\end{equation}
for any  $y\in V_m$ and $\nu\in m+{\Bbb N}-1.$

If, moreover, $k^\ast=0,$ then there exist $A>0,\,m\in{\Bbb N}$ such that
\begin{equation}\label{eq:xza1}
\exp(-A(ln(\nu)+\nu^{1-1/k}))(\rho_s)^{\nu}\omega_{n,y}(m)
\le
\omega_{n,y}(\nu)\end{equation}
for any  $y\in V_m$ and $\nu\in m+{\Bbb N}-1.$

\bf Proof. \rm Since $V_m=V_{m,1}^\wedge,$ it follows that
 the inequality (\ref{eq:za2}) holds with $\theta=1$ for any $y\in V_m.$
 For the full proof of the Lemma let us make some not large
 changes in the proof of the Lemma 2 in [\ref{r:ci0}].

The condition $k^\ast=0$ implies the inequality
\begin{equation}\label{eq:e}
a_0^\sim=T(0)\ne0,
\end{equation}
and for $p$ in Theorem 6 of the paper [\ref{r:cf0}] the equality $p=n.$ Let
\begin{equation}\label{eq:f}
A(\nu)=\left(\matrix
0&1&0&\ldots&0\\
0&0&1&\ldots&0\\
\vdots&\vdots&\vdots&\ldots&\vdots\\
0&0&0&\ldots&1\\
-a_0(\nu)&-a_1(\nu)&-a_2(\nu)&\ldots&-a_{n-1}(\nu)
\endmatrix\right),
\end{equation}
 where $\nu\in{\Bbb N}-1,$
 and let
\begin{equation}\label{eq:g}
A^\sim=\left(\matrix
0&1&0&\ldots&0\\
0&0&1&\ldots&0\\
\vdots&\vdots&\vdots&\ldots&\vdots\\
0&0&0&\ldots&1\\
-a_0^\sim&-a_1^\sim&-a_2^\sim&\ldots&-a_{n-1}^\sim.
\endmatrix\right)
\end{equation}
Let $\lambda_j,$ where $j=1,\,\ldots,\,r,$ is the sequence of all the
 mutually distinct roots of the polynomial (\ref{eq:zd})
 and $k^\ast_j$ is the multiplicity of the root $\lambda_j.$ Then
$$k=\sup\{k^\ast_j\colon j=1,\,\ldots,\,r\}.$$
 Clearly,
$$\sum\limits_{j=0}^r k^\ast_j=n,$$
$$
\rho_s\le\vert\lambda_j\vert\le\rho_1=\Vert A\Vert_{sp}<R=h^\sim(A^\sim)+1,
$$
 where $j=1,\,\ldots,\,r.$ In view of (\ref{eq:zc}),
 there exists $C_1>0$ such that
\begin{equation}\label{eq:h}
h^\sim(A(\nu)-A^\sim)\le C_1/(\nu+1).
\end{equation}
Therefore, according to the Lemma 2,
 for any $\varepsilon>0$ the following inequalities holds
\begin{equation}\label{eq:i}
p^\sim_{A^\sim,\varepsilon}(A(\nu)-A^\sim)\le
\end{equation}
$$(\gamma^\ast(A^\sim))^2(\max(\varepsilon,1/\varepsilon))^{k-1}
h^\sim(A(\nu)-A^\sim)\le$$
$$(\gamma^\ast(A^\sim))^2
(\max(\varepsilon,1/\varepsilon))^{k-1} C_1/(\nu+1),
$$
\begin{equation}\label{eq:aj}
p^\sim_{A^\sim,\varepsilon}(A(\nu))\le\rho_1+\sign(k-1)\varepsilon+
\end{equation}
$$(\gamma^\ast(A^\sim))^2
(\max(\varepsilon,1/\varepsilon))^{k-1}
C_1/(\nu+1)\le
$$
$$(\rho_1+(\sign(k-1))\varepsilon)\left(1+(\gamma^\ast(A^\sim))^2
(\max(\varepsilon,1/\varepsilon))^{k-1}\frac{C_1}{\rho_1(\nu+1)}\right),$$
\begin{equation}\label{eq:aa}
\ln(p^\sim_{A^\sim,\varepsilon}(A(\nu)))\le\ln(\rho_1)+
\ln(1+(\sign(k-1))C_{1,1}\varepsilon)+
\end{equation}
$$\ln\left(1+
(\max(\varepsilon,1/\varepsilon))^{k-1}C_{1,2}/(\nu+1)\right),$$
where
$$C_{1,1}=1/\rho_1,\,C_{1,2}=(\gamma^\ast(A^\sim))^2C_1/\rho_1$$
and $\nu\in{\Bbb N}.$

We consider first the case $k>1.$ For given $\nu\in{\Bbb N}+1$
 we take $d\in\Bbb N$ and $a_i\in[1,\nu]\cap\Bbb N,$
 where $i=0,\,\ldots,\,d,$ in such a way that
\begin{equation}\label{eq:ab}
a_0=1,\,a_d=\nu,a_{i-1}<a_i\le2a_{i-1},
\end{equation}
where $i=1,\,\ldots,\,d,$ and
\begin{equation}\label{eq:ac}
d\le\frac{\ln(\kappa)}{\ln(2)}+1.
\end{equation}
 According to the Corollary of the Lemma 1
 and (\ref{eq:ab}) -- (\ref{eq:ac}),
\begin{equation}\label{eq:ad}
\sum\limits_{\kappa=a_{i-1}}^{a_{i}-1}
\ln(p^\sim_{A^\sim,\varepsilon}(A(\kappa)))\le
\end{equation}
$$(a_i-a_{i-1})\ln(\rho_1)+(a_i-a_{i-1})\ln(1+C_{1,1}\varepsilon)\,+$$
$$\sum\limits_{\kappa=a_{i-1}}^{a_{i}-1}\ln\left(1+
C_{1,2}(\max(\varepsilon,1/\varepsilon))^{k-1}/(\kappa+1)\right)\le$$
$$(a_i-a_{i-1})(\ln(\rho_1)+\ln(1+\varepsilon/\rho_1))+
3C_{1,2}
(\max(\varepsilon,1/\varepsilon))^{k-1},$$
where $i=1,\,\ldots,\,d.$ We take now
 in (\ref{eq:ad}) $\varepsilon=\varepsilon_i=(a_{i-1})^{-1/k}.$
 Then we obtain the inequality
$$\sum\limits_{\kappa=a_{i-1}}^{a_{i}-1}
\ln(p^\sim_{A^\sim,\varepsilon}(A(\kappa)))\le$$
$$(a_i-a_{i-1})\ln(\rho_1)+(C_{1,1}+3C_{1,2})(a_{i-1})^{1-1/k}=$$
$$
(a_i-a_{i-1})\ln(\rho_1)+O(2^{(i-1)(1-1/k)}),$$
where $i=1,\,\ldots,\,d.$ Therefore
\begin{equation}\label{eq:ae}
\ln(p^\sim_{A^\sim,\varepsilon_i}\left(
\prod\limits_{\kappa=1}^{a_i-a_{i-1}}A(a_i-\kappa)\right)\le
\end{equation}
$$(a_i-a_{i-1})\ln(\rho_1)+O(2^{(i-1)(1-1/k)}$$
and, in view of (\ref{eq:ef}),
\begin{equation}\label{eq:af}
\ln(h^\sim\left(
\prod\limits_{\kappa=1}^{a_i-a_{i-1}}A(a_i-\kappa)\right)\le
\end{equation}
$$ \ln((\gamma^\ast(A^\sim))^2)(\max(\varepsilon,1/\varepsilon))^{k-1})+
(a_i-a_{i-1})\ln(\rho_1)+O(2^{(i-1)(1-1/k)})\le$$
$$
2\ln(\gamma^\ast(A^\sim))+(i-1)(1-1/k)ln(2)+
(a_i-a_{i-1})\ln(\rho_1)+O(2^{(i-1)(1-1/k)})=
$$
$$
(i-1)(1-1/k)\ln(2)+
(a_i-a_{i-1})\ln(\rho_1)+O(2^{(i-1)(1-1/k)}),
$$
where $i=1,\,\ldots,\,d.$ Consequently,
\begin{equation}\label{eq:ag}
\ln\left(h^\sim\left(
\prod\limits_{\kappa=1}^{nu-1} A(\nu-\kappa)\right)\right)=
\end{equation}
$$\ln\left(h^\sim\left(\left(\prod\limits_{i=1}^d
\prod\limits_{\kappa=1}^{a_i-a_{i-1}}A(a_i-\kappa)\right)\right)\right)=
$$
$$=\nu\ln(\rho_1)+O(\nu^{(1-1/k)}),$$
where $\nu\in{\Bbb N}.$

If $k=1,$ then, according to (\ref{eq:aa}),
$$\ln(p^\sim_{A^\sim,1}(A(\nu)))\le\ln(\rho_1)+C_{1,2}/(\nu+1),$$
\begin{equation}\label{eq:ah}
\ln
\left(h^\sim\left(\prod\limits_{\kappa=1}^{\nu-1}
A(\nu-\kappa)\right)\right)\le
\end{equation}
$$\ln\left((\gamma^\ast(A^\sim))^2p^\sim_{A^\sim,1}
\left(\prod\limits_{\kappa=1}^{\nu-1}A(\nu-\kappa)\right)\right)\le$$
$$\nu\ln(\rho_1)+O(\ln(e\nu)).$$
where $\nu\in{\Bbb N}.$ In view of (\ref{eq:ag}) -- (\ref{eq:ah}),
\begin{equation}\label{eq:ai}
\ln\left(h^\sim\left(
\prod\limits_{\kappa=1}^{\nu-1} A(\nu-\kappa)\right)\right)=
\end{equation}
$$=\nu\ln(\rho_1)+O((\nu+1)^{(1-1/k)})+O(\ln(e\nu))),$$
where $\nu\in{\Bbb N}.$

As in Section 3 of [\ref{r:cf0}],
 let $K$ denotes one of the fields $\Bbb R$ or $\Bbb C$ and $L$
 denotes a linear normed space over $K$ with norm $p=p(x).$ If $L=K^n,$
 we fix as $p=p(x),$ wehre $x\in K^n,$ the maximum of the
 absolute values of coordinates of the element $x$
 in the standard basis, i.e.
\begin {equation} \label{eq:bj}
p(x)=h(x)=\sup(\{\vert x_1\vert,\,\ldots,\,\vert x_n\vert\}),
\end{equation}
 where
$$x=\left(\matrix x_1\\\vdots\\ x_n\endmatrix\right).$$
If $L$ is a Banach space with the norm $p,$ then
$K-$algebra of all the linear continuous operators acting in $L$
 will be denoted by ${\mathfrak M^\wedge}(L),$ and the norm
 on ${\mathfrak M^\wedge}(L),$
 associated with the norm $p$ will be denoted by $p^\sim.$ So,
$$p^\sim(A)=sup(\{p(AX)\colon X\in L,\,p(X)\le1\}).$$
It is well known that the associciated with $h$ norm on $\Mat_n({\Bbb C})$
is defined as follows
\begin {equation} \label{eq:ba}
h^\sim(A)=sup\left(\left\{\sum\limits_{k=1}^n\vert a_{i,j}\vert
\colon i=1,\,\ldots,\,n\right\}\right),
\end{equation}
where $A=(a_{i,k})\in\Mat_n({\Bbb C}).$ The norms $h$ and $h^\sim$ coincide
 respectiely with with the norms $q_\infty$ and $q_\infty^\sim$
 considered in section 6 of the paper [\ref{r:cb}].

Let $m\in{\Bbb N},$ and let $E_m(L)$ be the set $L^{m-1+\Bbb N}$
 of all the maps of the set $m-1+\mathbb N$ into $L.$ The set $E_m(L)$
 is a linear space over $K,$ where the muliplication
 of the elements by the number from $K$ and addition
 of the elements are defined coordinate-wise. The subspace  of $E_m(L)$
 composed by all the constant maps is isomorphic to the space $L,$
 and we identify this subspace with $L.$

As in Section 4 of [\ref{r:cf0}],
 for any $y\in E_m({\Bbb C})$ and $n\in\Bbb N$ let $Y_{n,y}$
 and $Y_{n,y}^{\#}$ denote the elements in the space $E_m({\Bbb C}^n),$
 wich are determined respectively
 by means the following equalities:
\begin{equation}\label{eq:bb}
Y_{n,y}(\nu)=\left(\matrix y(\nu)\\\vdots\\ y(\nu+n-1)\endmatrix\right),
\end{equation}
\begin{equation}\label{eq:bc}
Y_{n,y}^{\#}(\nu)=\left(\matrix 0\\\vdots\\ 0\\ y(\nu)\endmatrix\right),
\end{equation}
where $\nu\in m-1+{\Bbb N}.$ Let is fixed $m\in{\Bbb N}.$
 If $y=y(\nu)$ is a solution of the equation (\ref{eq:zb})
 for $\nu\in {\Bbb N}+m-1,$ then
$$Y_{n,y}(\nu)=\left(\prod\limits_{\kappa=1}^{\nu-m} A(\nu-\kappa)\right)
Y_{n,y}(m),$$
where $\nu\in{\Bbb N},$ and, in view of (\ref{eq:ai}),
\begin{equation}\label{eq:bd}
\omega_{n,y}(\nu)=h\left(Y_{n,y}(\nu)\right)\le
\end{equation}
$$h^\sim\left(\prod\limits_{\kappa=1}^{\nu-m} A(\nu-\kappa)\right)
h\left(Y_{n,y}(m)\right)=
$$
$$\exp\left(O(1)\left(\nu^{1-1/k}+\ln(\nu)\right)\right)
\left(\rho_1\right)^\nu\omega_{n,y}(m),$$
where $\nu\in{\Bbb N}+m-1.$ So, with $m=1$ the asserted by the Lemma
 the upper estimate (\ref{eq:xza2}) of the value
 $\omega_{n,y}(\nu)\left(\rho_\theta\right)^{-\nu}=
\omega_{n,y}(\nu)\left(\rho_1\right)^{-\nu}$
 is obtained. We shall take now
\begin {equation} \label{eq:ei0}
\varepsilon\in(0,\rho_s/2).
\end{equation}
 Then, in view of (\ref{eq:ei}),
 \begin {equation} \label{eq:ei1}
 1/\rho_s\le(p_{A,\varepsilon})^\sim(A^{-1})
\le \end{equation}
$$1/\rho_s+
 2(\sign(k-1))\varepsilon/\rho_s^2\le2/\rho_s.
$$
Let is fixed
 \begin {equation} \label{eq:ei2}
m\in\max([(2/\rho_s)^k],\,
[2C_1h^\sim\left(\left(A^\sim\right)^{-1}\right)])+{\Bbb N}.
\end{equation}
Then
\begin {equation} \label{eq:ei4}
h^\sim\left(\left(A^\sim\right)^{-1}(A(\nu)-A^\sim)\right)
\le h^\sim\left(\left(A^\sim\right)^{-1}\right)C_1/(\nu+1)\le
\end{equation}
$$h^\sim\left(\left(A^\sim\right)^{-1}\right)C_1/
([2C_1h^\sim\left((A^\sim)^{-1}\right)]+1)\le1/2,$$
if $\nu\in m-1+{\Bbb N},$ the matrix
$E+\left((A^\sim)^{-1}(A(\nu)-A^\sim)\right)$ is invertible,
$$h^\sim\left(\left(E+\left((A^\sim)^{-1}(A(\nu)-
A^\sim)\right)\right)^{-1}\right)\le2,$$
if $\nu\in m-1+{\Bbb N},$
there exists the matix
$$\left(A(\nu)\right)^{-1}=
\left(E+\left(\left(A^\sim\right)^{-1}(A(\nu)-
A^\sim)\right)\right)^{-1}\left(A^\sim\right)^{-1},$$
if $\nu\in m-1+{\Bbb N},$
moreover,
$$h^\sim\left(\left(A(\nu)\right)^{-1}\right)\le
h^\sim\left(\left(E+\left((A^\sim)^{-1}(A(\nu)-
A^\sim)\right)\right)^{-1}\right)
h^\sim\left(\left(A^\sim\right)^{-1}\right)\le$$
$$
2h^\sim\left(\left(A^\sim\right)^{-1}\right)$$
and, finally,
$$h^\sim\left(\left(A(\nu)\right)^{-1}-
\left(A^\sim\right)^{-1}\right)=$$
$$h^\sim\left(\left(A(\nu)\right)^{-1}
\left(A^\sim-A(\nu)\right)\left(A^\sim\right)^{-1}\right)\le$$
$$h^\sim\left(\left(A(\nu)\right)^{-1}\right)
h^\sim\left(\left(A^\sim-A(\nu)\right)\right)
h^\sim\left(\left(A^\sim\right)^{-1}\right)\le C_2/(\nu+1),$$
where
$$C_2=\left(h^\sim\left(\left(A^\sim\right)^{-1}\right)\right)^2C_1$$
and $\nu\in m-1+{\Bbb N}.$
Therefore
\begin{equation}\label{eq:i1}
p^\sim_{A^\sim,\varepsilon}((A(\nu))^{-1}-(A^\sim)^{-1})\le
\end{equation}
$$(\gamma^\ast(A^\sim))^2(\max(\varepsilon,1/\varepsilon))^{k-1}
h^\sim((A(\nu))^{-1}-(A^\sim))^{-1}\le$$
$$
(\gamma^\ast(A^\sim))^2(\max(\varepsilon,1/\varepsilon))^{k-1}C_2/(\nu+1),
$$
where $\nu\in m-1+{\Bbb N}.$ In view of (\ref{eq:i1}) and (\ref{eq:ei}),
\begin{equation}\label{eq:aj1}
p^\sim_{A^\sim,\varepsilon}((A(\nu))^{-1})\le
1/\rho_s+\end{equation}
$$ 2(\sign(k-1))\varepsilon/\rho_s^2\le2/\rho_s+
(\gamma^\ast(A^\sim))^2
(\max(\varepsilon,1/\varepsilon))^{k-1}
C_2/(\nu+1)\le
$$
$$(1/\rho_s+2(\sign(k-1))\varepsilon/\rho_s^2)
\left(1+(\gamma^\ast(A^\sim))^2
(\max(\varepsilon,1/\varepsilon))^{k-1}\rho_sC_2\rho_s/(\nu+1)\right),$$
\begin{equation}\label{eq:aa2}
\ln(p^\sim_{A^\sim,\varepsilon}((A(\nu))^{-1}))\le\ln(1/\rho_s)+
\end{equation}
$$\ln(1+(\sign(k-1))C_{2,1}\varepsilon)+
\ln\left(1+
(\max(\varepsilon,1/\varepsilon))^{k-1}C_{2,2}/(\nu+1)\right),$$
where
$$C_{2,1}=2/\rho_s,\,C_{2,2}=(\gamma^\ast(A^\sim))^2C_2\rho_s$$
and $\nu\in{\Bbb N}-1+m.$ We take $\nu\in m-1+{\Bbb N}.$

We consider first the case $k>1$ again now. For given $\nu\in{\Bbb N}+m$
 we take $d\in\Bbb N$ and $a_i\in[m,\nu]\cap\Bbb N,$
where $i=0,\,\ldots,\,d,$ in such a way that
\begin{equation}\label{eq:ab1}
a_0=m,\,a_d=\nu,a_{i-1}<a_i\le2a_{i-1}\le m2^i,
\end{equation}
where $i=1,\,\ldots,\,d,$ and
\begin{equation}\label{eq:ac1}
d\le\frac{\ln(\nu)}{\ln(2)}+1.
\end{equation}
According to the Corollary of the Lemma 1
 and (\ref{eq:ab1}) -- (\ref{eq:ac1}),
\begin{equation}\label{eq:ad1}
\sum\limits_{\kappa=a_{i-1}}^{a_{i}-1}
\ln(p^\sim_{A^\sim,\varepsilon}((A(\kappa))^{-1}))\le
\end{equation}
$$(a_i-a_{i-1})\ln(1/\rho_s)+(a_i-a_{i-1})\ln(1+C_{2,1}\varepsilon)\,+$$
$$\sum\limits_{\kappa=a_{i-1}}^{a_{i}-1}
\ln\left(1+(\max(\varepsilon,1/\varepsilon))^{k-1}
C_{2,2}/(\kappa+1)\right)\le$$
$$(a_i-a_{i-1})(\ln(1/\rho_s)+\ln(1+C_{2,1}\varepsilon))+
3(\gamma^\ast(A^\sim))^2
(\max(\varepsilon,1/\varepsilon))^{k-1}C_{2,2},$$
where $i=1,\,\ldots,\,d.$ We take now in (\ref{eq:ad1})
 $\varepsilon=\varepsilon_i=(a_{i-1})^{-1/k}.$
 Then, in view of (\ref{eq:ei2}),
$$a_{i-1}\ge m>(2/\rho_s)^k,\,(a_{i-1})^{1/k}>2/\rho_s,\,
\varepsilon_i<\rho_s/2,$$
where $i=1,\,\ldots,\,d,$ and therefore (\ref{eq:ei0})
 and (\ref{eq:ad1}) hold. Consequently,
\begin{equation}\label{eq:ad2}
\sum\limits_{\kappa=a_{i-1}}^{a_{i}-1}
\ln(p^\sim_{A^\sim,\varepsilon}((A(\kappa))^{-1}))\le
\end{equation}
$$(a_i-a_{i-1})\ln(1/\rho_s)+(C_{2,1}+3C_{2,2})(a_{i-1})^{1-1/k}=$$
$$
(a_i-a_{i-1})\ln(1/\rho_s)+O(2^{(i-1)(1-1/k)}),$$
where $i=1,\,\ldots,\,d.$ Therefore
\begin{equation}\label{eq:ae1}
\ln(p^\sim_{A^\sim,\varepsilon_i}\left(
\prod\limits_{\kappa=a_{i-1}^{a_i-1}}(A(\kappa))^{-1}\right)\le
\end{equation}
$$(a_i-a_{i-1})\ln(\rho_1)+O(2^{(i-1)(1-1/k)}$$
and,in view of (\ref{eq:ef}),
\begin{equation}\label{eq:af1}
\ln(h^\sim\left(
\prod\limits_{\kappa=1}^{a_i-a_{i-1}}(A(\kappa))^{-11}\right)\le
\end{equation}
$$ \ln((\gamma^\ast(A^\sim))^2)(\max(\varepsilon_i,1/\varepsilon_i))^{k-1})+
(a_i-a_{i-1})\ln(\rho_1)+O(2^{(i-1)(1-1/k)})\le$$
$$
2\ln(\gamma^\ast(A^\sim))+(i-1)(1-1/k)ln(2)+
(a_i-a_{i-1})\ln(\rho_1)+O(2^{(i-1)(1-1/k)})=
$$
$$
(i-1)(1-1/k)ln(2)+
(a_i-a_{i-1})\ln(\rho_1)+O(2^{(i-1)(1-1/k)}),
$$
where $i=1,\,\ldots,\,d.$ Consequently,
\begin{equation}\label{eq:ag1}
\ln\left(h^\sim\left(
\prod\limits_{\kappa=m}^{\nu-1}(A(\kappa))^{-1}\right)\right)=
\end{equation}
$$\ln\left(h^\sim\left(\left(\prod\limits_{i=1}^d
\prod\limits_{\kappa=1}^{a_i-a_{i-1}}(A(\kappa))^{-1}\right)\right)\right)=
$$
$$=\nu\ln(1/\rho_s)+O(\nu^{(1-1/k)}),$$
where $\nu\in{\Bbb N}+m.$

If $k=1,$ then, according to (\ref{eq:aa2}),
$$\ln(p^\sim_{A^\sim,1}(A(\nu)))\le\ln(1\rho_s)+C_{2,2}/(\nu+1),$$
\begin{equation}\label{eq:ah1}
\ln
\left(h^\sim\left(\prod\limits_{\kappa=m}^{\nu-1}
(A(\kappa))^{-1}\right)\right)\le
\end{equation}
$$\ln\left((\gamma^\ast(A^\sim))^2p^\sim_{A^\sim,1}
\left(\prod\limits_{\kappa=m}^{\nu-1}(A(\kappa))^{-1}\right)\right)\le$$
$$\nu\ln(1/\rho_s)+O(\ln(e\nu)).$$
where $\nu\in{\Bbb N}.$ In view of (\ref{eq:ag1}) -- (\ref{eq:ah1}),
\begin{equation}\label{eq:ai1}
\ln\left(h^\sim\left(
\prod\limits_{\kappa=m}^{\nu-1} A(\nu-\kappa)\right)\right)=
\end{equation}
$$=\nu\ln(1/\rho_s)+O(\nu^{(1-1/k)})+O(\ln(e\nu))),$$
where $\nu\in{\Bbb N}+m.$ Since
$$Y_{n,y}(m)=
\left(\prod\limits_{\kappa=1}^{\nu-m} A(\nu-\kappa)\right)^{-1}
Y_{n,y}(\nu)=$$
$$\left(\prod\limits_{\kappa=m}^{\nu-1}\left(A(\kappa)\right)^{-1}\right)
Y_{n,y}(\nu),$$
it follows that
\begin{equation}\label{eq:bf1}
\omega_{n,y}(m)=h\left(Y_{n,y}(m)\right)\le
\end{equation}
$$h^\sim\left(\prod\limits_{\kappa=m}^{\nu-1}
\left(A(\kappa)\right)^{-1}\right)
h\left(Y_{n,y}(\nu)\right)=
$$
$$\exp\left(O(1)\left(\nu^{1-1/k}+\ln(\nu)\right)\right)
\left(\rho_1\right)^{-\nu}\omega_{n,y}(\nu),$$
where $\nu\in{\Bbb N}+m.$ So, with $m$ from  (\ref{eq:ei2})
 the asserted by the Lemma the lower estimate of the value
 $\omega_{n,y}(\nu)\left(\rho_s\right)^{-\nu}$ is obtained. $\blacksquare$

\bf Remark 1. \it  In the case $s=1,\,k^\ast=0$ the assertions of
 the Lemma and Theorem 7 coincide. \rm

\bf Lemma 4 (Perron's decomposition Lemma, [\ref{r:a}], Hilfsatz 3)). \it
 Let the characteristic polynomial (\ref{eq:zd}) is represented as product
\begin{equation}\label{eq:1h}
T(z)=T_1(z)T_2(z),
\end{equation}
where
\begin{equation}\label{eq:1i}
T_1(z)=\sum\limits_{\alpha=0}^pb_\alpha^\sim z^{\alpha},\,
T_2(z)=\sum\limits_{\beta=0}^qu_\beta^\sim z^{\beta},\end{equation}
with $b_p^\sim=u^\sim_q=a^\sim_n=1$ and absolute value of each root
 of $T_1(z)$ is greater than the absolute value of each root of $T_2(z).$

Then there exist $m\in{\Bbb N},\,
b_\alpha(\nu)\in{\Bbb C},\alpha=0,\,\ldots,\,p,\nu\in{\Bbb N}+m-1,$\\
and
$u_\beta(\nu)\in{\Bbb C},\beta=0,\,\ldots,\,q,\nu\in{\Bbb N}+m-1,$
such that
\begin{equation}\label{eq:1aa}\lim\limits_{\nu\to\infty}b_\alpha(\nu)=
b_\alpha^\sim,\,\alpha=0,\,\ldots,\,p,\,b_p(\nu)=1,b_0(\nu)\ne0,
\end{equation}
\begin{equation}\label{eq:1ab}\lim\limits_{\nu\to\infty}u_\beta(\nu)=
u_\beta^\sim,\beta=0,\,\ldots,\,q,\,u_q(\nu)=1,\end{equation}
where $\nu\in{\mathbb N}+m-1,$ and, moreover, the equation (\ref{eq:zb})
 is equivalent to the equation
$$\sum\limits_{\alpha=0}^pb_\alpha(\nu)y(\nu+\alpha)=r(\nu),$$
where $\nu\in{\mathbb N}-1+m$ and $r(\nu)$ satisfies to the equation
$$\sum\limits_{\beta=0}^qu_\beta(\nu)r(\nu+beta)=0$$
with $\nu\in{\mathbb N}-1+m.$

\bf Lemma 5. \it
Let all the conditions of the Perron's decomposition Lemma are fulfilled and
$a_k(\nu)-a_k^\sim=O(1/(\nu+1)),k=0,\,\ldots,\,n,$
when $\nu\to\infty.$

Then for $b_\alpha(\nu)$ with $\alpha=0,\,\ldots,\,p,$ and
$u_\beta(\nu)$ with $\beta=0,\,\ldots,\,q,$ from the assertion of
the Perron's decomposition Lemma the following conditions
 are fulfilled (cf. with (\ref{eq:1aa}) and (\ref{eq:1ab})):
$$ b_\alpha(\nu)-b_\alpha^\sim=O(1/\nu),\,
\alpha=0,\,\ldots,\,p,b_0(\nu)=1,\,b_q(\nu)\ne0,$$
$$u_\beta(\nu)-u_\beta^\sim=O(1/\nu),\,
\beta=0,\,\ldots,\,q,\,u_0(\nu)=1,$$
where $\nu\in{\mathbb N}-1+m.$

\bf Proof. \rm The Lemma is direct corollary of
the Theorem 5 in [\ref{r:cb}]. $\blacksquare$
{\begin{center}\large\bf \S 2.
 Proof of the theorem 8.\end{center}}
We use below the notations of the section 3 in [\ref{r:ci0}].
 Let $K$ be one of the fields $\Bbb R$ or $\Bbb C.$ Let $L$
 be a linear normed space over  $K$ with a norm $p(x).$ In the case $L=K^n$
 we fix as $p(x),$ wehre $x\in K^n,$ the maximum of the absolute values
 ofcoordinates of $x$ in the standard basis, i.e.
\begin {equation}
\label{eq:2cea}
p(x)=h(x)=\sup(\{\vert x_1\vert,\,\ldots,\,\vert x_n\vert\}),
\end{equation}
 where
$$x=\left(\matrix x_1\\\vdots\\ x_n\endmatrix\right).$$
If $L$ is a Banach space with the norm $p,$ then
 $K-$algebra of all the linear continuous operators acting in $L$
 will be denoted by ${\mathfrak M^\wedge}(L),$ and the norm
 on ${\mathfrak M^\wedge}(L),$ associated with the norm $p$
 will be denoted by $p^\sim.$ So,
$$p^\sim(A)=sup(\{p(AX)\colon X\in L,\,p(X)\le1\}).$$
It is well known that the associciated with $h$ norm on $\Mat_n({\Bbb C})$
 is defined as follows
\begin {equation} \label{eq:2ceb}
h^\sim(A)=sup\left(\left\{\sum\limits_{k=1}^n\vert a_{i,j}\vert
\colon i=1,\,\ldots,\,n\right\}\right),
\end{equation}
where $A=(a_{i,k})\in\Mat_n({\Bbb C}).$
 The norms $h$ and $h^\sim$ coincide respectiely with
 the norms $q_\infty$ and $q_\infty^\sim$ considered
 in [\ref{r:cb}], section 6. Let $m\in{\Bbb N}-1,$ and we denote by $E_m(L)$
 the set $L^{m-1+\Bbb N}$ of all the
 maps of the set $m-1+\mathbb N$ into $L.$
 The set $E_m(L)$ is a linear space over $K,$
 where the muliplication of the elements by the number from $K$ and addition
 of the elements is defined coordinate-wise. The subspace  of $E_m(L)$
 composed by all the constant maps is isomorphic to $L,$ and we identify
 this subspace with $L.$

We denote by ${\mathfrak M}^\vee(L)$ the space of all the $K-$linear
 maps of the space $L$ in $L.$ If $\phi\in{\mathfrak M}^\vee(L)$
 and $\psi\in{\mathfrak M}^\vee(L),$ then $\phi\circ\psi$ denotes
 the composition of operators $\phi$ and $\psi,$ so that
$(\phi\circ\psi)f=\phi((\psi f))$ for each $f\in L.$
 For $x\in E_m(L)$ let
$$p_{m,\infty}(x)=\sup(\{p(x(\nu))\colon\nu\in m-1+\Bbb N\}.$$
Let further
$$E_{m,\infty}(L)=\{x\in E_m(L)\colon p_{m,\infty}(x)\ne\infty \},$$
$$E_{m,0}(L)=\{x\in E_m(L)\colon \lim_{\nu\to\infty} p(x(\nu))=0\},$$
$$E_m^\to(L)=L+E_{m,0}(L).$$
Clearly, the space $E_m^\to(L)$ consists of all the $y\in E_m(L),$ for which
 there exists
$$\lim(y)=\lim\limits_{\nu\to\infty}(y(\nu)).$$
Let $m\in{\mathbb N}-1,\,\mu\in m-1+\Bbb N$ and ler $r_{m,\mu}$
 be the operator of restriction of the elements $y\in E_m(L)$
 on te set $m-1+\Bbb N.$ Clearly, the map $r_{m,\mu}$ is an epimorphism
 of the space $E_m(L)$ onto the space $E_\mu(L).$ If $L$ is a $K$-algebra,
 then $E_m(L)$ is a $K$-algebra, where the muliplication
 and addition of the elements is defined coordinate-wise; so,
 in this case $r_{m,\mu}$ is an epimorphism
 of $K$-algebra $E_m(L)$ onto $K$-algebra $E_\mu(L).$
 If $L$ be an algebra with unity, let $L^\ast$ denotes the group of
 all its invertible elements. Then
$$(L^\ast)^{m-1+\Bbb N}\subset L^{m-1+\Bbb N};$$
 we denote below $(L^\ast)^{m-1+\Bbb N}$ by $E_m(L^\ast).$
 Clearly,
$$E_m(L^\ast)=(E_m(L))^\ast.$$
Let $L={\mathbb C}^n,\,y\in E_m(L),$ and let $y_i(\nu)$
 denotes the $i$-th coordinate of the element $y(\nu),$ where
 $i=1,\,\ldots,\,n,\,\nu\in m-1+\Bbb N;$ then the space $(E_m({\Bbb C}))^n$
 contains an element $\omega(y),$ which has $y_i(\nu)$ as the value of
 its $i$-th coordinate at the point $\nu\in m-1+\Bbb N.$
  So we obtain the natural isomorphism $\omega$
 of the algebra $E_m({\Bbb C}^n)$ onto $(E_m({\Bbb C}))^n.$
 This map $\omega$ induces an isomorphism
 of the algebra $E_m(\Mat_n({\Bbb C}))$ onto $\Mat_n(E_m({\Bbb C})).$
 If $L$ is a $K-algebra,$ then
 each element $a\in E_m(L)$ determines an acting
 on the space $E_m(L)$ $K$-linear
 operator $\mu_a\in{\mathfrak M}^\vee(E_m(L)),$
 which turns any $y\in E_m(L)$ into $\mu_ay=ay.$
 On the space $E_m(L)$ acts
 also $K$-linear operator $\bigtriangledown\in{\mathfrak M}^\vee(L),$ which
 turns any element $y\in E_m(L)$ in the $\bigtriangledown y\in E_m(L)$
 such that
$$(\bigtriangledown y)(\nu)=y(\nu+1)$$
 for any $\nu\in m-1+\Bbb N.$ Let us consider
 the subring ${\mathfrak A}_m(L)$ of the ring $\mathfrak M^\vee(L)$
 generated by the operator $\bigtriangledown$ and
 by all the operators $\mu_a,$ where $a\in E_m(L).$
 Clearly,
\begin{equation}\label{eq:2cbe}
\mu_a\circ\bigtriangledown^r\circ
\mu_b\circ\bigtriangledown^s=
\mu_{a\bigtriangledown^rb_k}\circ\bigtriangledown^{r+s},
\end{equation}
 where $\{r,s\}\subset {\Bbb N}-1,\,\{a,\,b\}\subset E_m(L).$
 For each $\alpha\in{\mathfrak A}_m(L)\diagdown\{0_{{\mathfrak A}_m(L)}\}$
 are uniquelly defined the number $\deg(\alpha)$ and representation
 of $\alpha$ in the form
\begin{equation}\label{eq:2cbf}
\alpha=\sum\limits^{\deg(\alpha)}_{k=0}\mu_{a_k}\circ
\bigtriangledown^k,
\end{equation}
where $a_k\in E_m(L)$ for $k=0,\,\ldots,\,\deg(\alpha)$
 and $a_{\deg(\alpha)}\ne0_{E_m(L)}.$
 Clearly, (\ref{eq:2cbf}) may be rewritten in the form
\begin{equation}\label{eq:2cbg}
\alpha=\sum\limits^{\infty}_{k=0}\mu_{a_k}\circ
\bigtriangledown^k,
\end{equation}
where $a_k=0_{E_m(L)}$ for $k\in\deg(\alpha)+\mathbb N.$
 It follows from (\ref{eq:2cbe}) that
 ${\mathfrak A}_m(L)$ is a graduated algebra, and if
\begin{equation}\label{eq:2cbh}
\beta=\sum\limits^p_{r=0}\mu_{b_r}\circ
\bigtriangledown^r\in{\mathfrak A}_m(L),
\end{equation}
\begin{equation}\label{eq:2cbi}
\gamma=\sum\limits^q_{s=0}\mu_{c_s}\circ
\bigtriangledown^s\in{\mathfrak A}_m(L),
\end{equation}
then
\begin{equation}\label{eq:2ccj}
\beta\gamma=\sum\limits^{p+q}_{k=0}
\sum\Sb\le r\le p\\0\le s\le q\\r+s=k\endSb
\mu_{b_r\bigtriangledown^rc_s}\circ\bigtriangledown^{r+s};
\end{equation}
 clearly, $\deg(\beta\gamma)=\deg(\beta)+\deg(\gamma),$
 if $b_p(\nu)^rc_q(\nu+p)$ is different from $0$
 at least for one $\nu\in m-1+\Bbb N.$ Let ${\mathfrak A}^\to_m(L)$
 be the ring generated by the operator $\bigtriangledown$ and
 by all the operators $\mu_a,$ where $a\in E^\to_m(L).$
 Since $\bigtriangledown a\in E^\to_m(L),$ if $a\in E^\to_m(L),$
 it follows, in view of (\ref{eq:2cbe}), that
 ${\mathfrak A}^\to_m(L)$ is a graduated subalgebra  ${\mathfrak A}^\to_m(L)$
 of the algebra  ${\mathfrak A}_m(L),$
 each $\alpha\in{\mathfrak A}_m(L)\diagdown\{0_{{\mathfrak A}_m(L)}\}$
 admits a representation in the form (\ref{eq:2cbf}) with
 $a_k\in E^\to_m(L)$ for $k=0,\,\ldots,\,\deg(\alpha)$
 and $a_{\deg(\alpha)}\ne0_{E_m(L)};$ to each such $\alpha$ corresponds
 the limit operator
\begin{equation}\label{eq:2cca}
\lim(\alpha)=\sum\limits^{\deg(\alpha)}_{k=0}\mu_{\lim(a_k)}\circ
\bigtriangledown^k,
\end{equation}
and polynomial
\begin{equation}\label{eq:2ccb}
P(\alpha,z)=\sum\limits^{\deg(\alpha)}_{k=0}\lim(a_k)z^k\in L[z].
\end{equation}
If $\alpha=0_{{\mathfrak A}_m(L)},$ then we put
$$\lim(\alpha)=0_{{\mathfrak A}_m(L)},\,P(\alpha,z)=0_{L[z]}.$$
The equality (\ref{eq:2cbe}) shows that the map
\begin{equation}\label{eq:2ccc}
\alpha\to P(\alpha,z)
\end{equation}
is an epimorphism of the algebra ${\mathfrak A}^\to_m(L)$ on the algebra $L[z]$
 (if the algebra $L$ is noncommutative, then we can treat the algebra $L[z]$ as
 a semigroup ring of the semigroup $({\Bbb N}-1,\,+)$ over the algebra $L$).
 We note that, if $\alpha\in {\mathfrak A}_m(\mathbb C),$ then,
 clearly, $\Ker(\alpha)$ coincides with the linear space
 of all the solutions of the equation (\ref{eq:zb}),
 and, moreover, if $\alpha\in {\mathfrak A}^\to_m(\mathbb C),$
 then the corresponding to $\alpha$ equation (\ref{eq:zb})
 is an equation of the Poincar'e type and $P(\alpha,z)$
 is its characterictical polynomial.

Let $L$ be an algebra with unity. The set of all
 $\alpha\in{\mathfrak A}_m(L)\diagdown\{0_{{\mathfrak A}_m(L)}\},$
 which have the representation (\ref{eq:2cbf}) with
$a_{\deg(\alpha)}\in E_m(L^\ast)$
 will be denoted further by ${\mathfrak A}_m(L)^\circ.$
 The set of all the the elements
 $\alpha\in{\mathfrak A}_m(L)\diagdown\{0_{{\mathfrak A}_m(L)}\},$
 which have the representation (\ref{eq:2cbf})
 with $a_{\deg(\alpha)}=1_{E_m(L)}$
 will be denoted further by ${\mathfrak A}_m(L)^\vee.$
 The set of all the
 $\alpha\in{\mathfrak A}_m(L)\diagdown\{0_{{\mathfrak A}_m(L)}\},$
 which have the representation (\ref{eq:2cbf}) with $a_0\in E_m(L^\ast),$
 will be denoted further by ${\mathfrak A}_m(L)^\wedge.$
 Let further
$$
{\mathfrak A}_m(L)^{\circ\wedge} ={\mathfrak A}_m(L)^\circ\bigcap
{\mathfrak A}_m(L)^\wedge,\,
{\mathfrak A}_m(L)^{\vee\wedge}={\mathfrak A}_m(L)^\vee\bigcap
{\mathfrak A}_m(L)^\wedge,
$$
$$
{\mathfrak A}_m(L)^{\circ\to} ={\mathfrak A}_m(L)^\circ\bigcap
{\mathfrak A}_m(L)^\to,\,
{\mathfrak A}_m(L)^{\vee\to}={\mathfrak A}_m(L)^\vee\bigcap
{\mathfrak A}_m(L)^\to,
$$
$$
{\mathfrak A}_m(L)^{\circ\wedge\to} ={\mathfrak A}_m(L)^{\circ\wedge}\bigcap
{\mathfrak A}_m(L)^\to,\,
{\mathfrak A}_m(L)^{\vee\wedge\to}={\mathfrak A}_m(L)^{\vee\wedge}\bigcap
{\mathfrak A}_m(L)^\to,
$$
Clearly, ${\mathfrak A}_m(L)^\circ$ consists of epimorphisms
 of the space $E_m(L)$ onto $E_m(L).$ The above map $r_{m,\mu}$
 induces epimorphism $r^\vartriangleright_{m,\mu}$ of
 the algebra ${\mathfrak A}_m(L)$ on the algebra ${\mathfrak A}_\mu(L)$ defined
 as follows:

if $\alpha\in{\mathfrak A}_m(L),$
\begin{equation}\label{eq:2ccd}
\alpha=\sum\limits^n_{k=0}\mu_{a_k}\circ
\bigtriangledown^k,
\end{equation}
then
\begin{equation}\label{eq:2cce}
r^\vartriangleright_{m,\mu}(\alpha)=
\sum\limits^n_{k=0}\mu_{r_{m,\mu}(a_k)}\circ
\bigtriangledown^k,
\end{equation}
where the operator $\bigtriangledown$ in (\ref{eq:2ccd}) acts in $E_m(L)$
 and the operator $\bigtriangledown$ in (\ref{eq:2cce}) acts in $E_\mu(L).$
 Clearly, $r^\vartriangleright_{m,\mu}$ surjectively maps

${\mathfrak A}_m(L)^\circ$ onto ${\mathfrak A}_\mu(L)^\circ,\,
{\mathfrak A}_m(L)^\vee$ onto ${\mathfrak A}_\mu(L)^\vee,$

${\mathfrak A}_m(L)^\wedge$ onto ${\mathfrak A}_mu(L)^\wedge,\,
{\mathfrak A}_m(L)^{\circ\wedge}$ onto ${\mathfrak A}_\mu(L)^{\circ\wedge},$

${\mathfrak A}_m(L)^{\vee\wedge}$ onto ${\mathfrak A}_\mu(L)^{\vee\wedge},$
${\mathfrak A}_m(L)^{\circ\to}$ onto ${\mathfrak A}_\mu(L)^{\circ\to},$

${\mathfrak A}_m(L)^{\vee\to}$ onto ${\mathfrak A}_\mu(L)^{\vee\to},$
${\mathfrak A}_m(L)^{\circ\wedge\to}$ onto
${\mathfrak A}_\mu(L)^{\circ\wedge\to},$

${\mathfrak A}_m(L)^{\vee\wedge\to}$ onto
${\mathfrak A}_\mu(L)^{\vee\wedge\to}.$
 Since the diagram
$$
\matrix
{E_m(L)}&{@>r_{m,\mu}>>}&
{E_\mu(L)}\quad\\
{\alpha\bigg\downarrow}&&
{\bigg\downarrow r^\vartriangleright_{m,\mu}(\alpha)}\\
{E_m(L)}&{@>>r_{m,\mu}>}&{E_\mu(L)}\quad\\
\endmatrix
$$
is commuative and therefore
\begin{equation}\label{eq:2ccf}
r_{m,\mu}\alpha=r^\vartriangleright_{m,\mu}(\alpha)r_{m,\mu},
\end{equation}
 it follows that $r_{m,\mu}$ surjectively maps $\Ker(r_{m,\mu}\alpha)$
 onto
$$\Ker(r^\vartriangleright_{m,\mu}(\alpha))\supset r_{m,\mu}\Ker(\alpha).$$

\bf Lemma 6. \it If $\mu\in m-1+\Bbb N$ and
 $\alpha\in{\mathfrak A}_m(L)^\wedge,$ then the operator $\alpha$
 bijectively maps $\Ker(r_{m,\mu})$ onto $\Ker(r_{m,\mu})$.

\bf Proof. \rm Proof is given in [\ref{r:cf0}], Lemma 3.

\bf Corollary 1. \it Let $\mu\in m-1+\Bbb N$
 and let $\alpha\in{\mathfrak A}_m(L)^\wedge.$
 If $$g\in E_m(L),\,x\in E_\mu(L),\,m\le\mu,\,
\alpha\in{\mathfrak A}_m(L)^\wedge,$$
$$r_{m,\mu}(g)=(r^\vartriangleright_{m,\mu}(\alpha))(x),$$
 then there exists a unique $y\in E_m(L)$ such that
$$\alpha(y)=g,\,r_{m,\mu}(y)=x;$$

\bf Proof. \rm Proof is given in [\ref{r:cf0}], Corollary 1 to the Lemma 3.

\bf Corollary 2. \it Let $\mu\in m-1+{\Bbb N}$
 and $\alpha\in{\mathfrak A}_m(L)^\wedge.$
 Then and $r_{m,\mu}$ bijectively maps $\Ker(\alpha)$ onto
$\Ker(r^\vartriangleright_{m,\mu}(\alpha))=r_{m,\mu}(\Ker(\alpha)).$

\bf Proof. \rm Proof is given in [\ref{r:cf0}], Corollary 2 to the Lemma 3.

If for the equation (\ref{eq:za0}) are fulfilled the conditions (\ref{eq:zc})
then
\begin{equation}\label{eq:2cci}
a_k=(a_k(0),\,a_k(1),\,\ldots,\,a_k(\nu),\,\ldots)\in E_0^\to(\mathbb C),
\end{equation}
where $k=0,\,\ldots,\,n.$ Moreover,
 $a_n=1_{E_0(\mathbb C)},$ for $\alpha$ in (\ref{eq:2ccd})
 $\Ker(\alpha)$ coincides with the linear over ${\mathbb C}$
 space of all the solutions of the equation (\ref{eq:za0}),
 polynomial (\ref{eq:zd}) is equal
 to the polynomial $P(\alpha,z)=P(r^\vartriangleright_{0,m}(\alpha),z),$
 where $m\in{\Bbb N},$ and the set $\Ker(r^\vartriangleright_{0,m}(\alpha))$
 coincides with linear over ${\mathbb C}$ space $V_m$ of all the solutions
 of the equation (\ref{eq:zb}).

 Let $\mathfrak v$ be the element in $E_{0,0},$ for which
$${\mathfrak v}(\nu)=\frac1{\nu+1},$$
 where $\nu\in{\Bbb N}-1$.
 Clearly,
 $r_{0,m}({\mathfrak v})E_{m,\infty}({\Bbb C})\subset E_{m,0}({\Bbb C})$
 for any $m\in{\Bbb N}-1.$ Let
$$E_{m,0}^\succ(L)=r_{0,m}({\mathfrak v})E_{m,\infty}(L),\,
E_m^\succ(L)=L+E_{m,0}^\succ(L).$$
Let us consider the ring ${\mathfrak A}_m^\succ(L)$ generated by the operator
 $\bigtriangledown$ and by all the operators $\mu_a,$ where
 $a\in E_m(L)^\succ.$ The Lemma 5 may be reformulated as follows:

\bf Lemma 7. \it Let
$\alpha\in{\mathfrak A}_0^\succ({\Bbb C})
\bigcap{\mathfrak A}_0^\vee({\Bbb C}),$ and $P(\alpha,z)$ coincides
 with the polynomial $T(z)$ in (\ref{eq:zd}) and (\ref{eq:1h}).

Then there exist $m\in\Bbb N$ and representation of the operator
 $r^\vartriangleright_{0,m}(\alpha)$ in the form
\begin{equation}\label{eq:2cdj}
r^\vartriangleright_{0,m}(\alpha)=\eta\beta
\end{equation}
such that
\begin{equation}\label{eq:2cda}
\eta\in{\mathfrak A}_m^\succ({\Bbb C})
\bigcap{\mathfrak A}_m^\vee({\Bbb C}),\,deg(\eta)=q
\end{equation}
\begin{equation}\label{eq:2cdb}
\beta\in{\mathfrak A}_m^\succ({\Bbb C})
\bigcap{\mathfrak A}_m^{\vee\wedge}({\Bbb C}),\,
\deg(\beta)=p=n-q,
\end{equation}
and the polynomials $P(\beta,z)$ , $P(\eta,z)$  coincide
 respectively with the polynomials $T_1(z)$ , $T_2(z)$ in (\ref{eq:1i}) .

\bf Lemma 8. \it Let are fulfilled all the conditions of the Theorem 7.
Let
$\alpha\in{\mathfrak A}_0^\succ({\Bbb C})
\bigcap{\mathfrak A}_0^\vee({\Bbb C})$ corresponds
 to the equation (\ref{eq:za0}), i.e. (with $m=0$) (\ref{eq:2ccd}) holds
\begin{equation}\label{eq:z2ccd}
\alpha=\sum\limits^n_{k=0}\mu_{a_k}\circ
\bigtriangledown^k,
\end{equation}
where
$a_n=1_{E_0(\mathbb C)}$ and
\begin{equation}\label{eq:z2cci}
a_k=(a_k(0),\,a_k(1),\,\ldots,\,a_k(\nu),\,\ldots)\in E_0^\succ(\mathbb C),
\end{equation}
for $k=0,\,\ldots,\,n.$

 Let the characteristic polynomial (\ref{eq:zd}) is represented as product
\begin{equation}\label{eq:2a}
P(\alpha,z)=T(z)=\prod\limits_{k=1}^{s+1}T_i(z),
\end{equation}
where
\begin{equation}\label{eq:2i}
T_i(z)=\sum\limits_{\alpha=0}^{e_i}b_{i,\alpha}^\sim z^{\alpha},\,
\end{equation}
with $b_{i,e_i}=a^\sim_n=1$ and absolute value of each root
 of the polynomial $T_i(z)$ is equal to $\rho_i.$

Then there exist $m\in\Bbb N$ and representation of the operator
$r^\vartriangleright_{0,m}(\alpha)$ in the form
\begin{equation}\label{eq:z2cdj}
r^\vartriangleright_{0,m}(\alpha)=
\prod\limits_{i=0}^s\beta_{s+1-i}
\end{equation}
such that
\begin{equation}\label{eq:z2cdb}
\beta_i\in{\mathfrak A}_m^\succ({\Bbb C})
\bigcap{\mathfrak A}_m^\vee({\Bbb C}),\,
\deg(\beta_i)=e_i,
\end{equation}
$$P(\beta_i,z)=T_i(z),$$
where $i=1,\,\ldots,\,s+1$
and
\begin{equation}\label{eq:z2cdc}
\beta_i\in{\mathfrak A}_m^\wedge({\Bbb C}),
\end{equation}
where $i=1,\,\ldots,\,s.$

\bf Proof. \rm The assertion of the Lemma may be obtained by means of the
 sequentielly applying of the Lemma 7. $\blacksquare$

Let
\begin{equation}\label{eq:xz2cdj}
\beta^\ast_\theta=
\prod\limits_{i=0}^{\theta-1}\beta_{\theta-i},\,
\beta^\wedge_\theta=
\prod\limits_{i=s+1-\theta}^{s}\beta_{s+1-i},
\end{equation}
where $\theta=1,\,\ldots\,s+1.$ In view of (ref{eq:2cdj}),
\begin{equation}\label{eq:yz2cdj}
r^\vartriangleright_{0,m}(\alpha)=
\beta^\wedge_{\theta+1}\beta^\ast_\theta,
\end{equation}
$\theta=1,\,\ldots\,s.$

 Let $C>0,\,n\in\Bbb N,$
\begin{equation}\label{eq:zcdd3}
w_{C,n}(\nu)=\left(\frac C{\nu}\right)^{\nu/n},
\end{equation}
where $\nu\in{\Bbb N};$
let further $\rho>0$
and
\begin{equation}\label{eq:zcdd4}
v_{C,n,\rho}(\nu)=
\rho^\nu\exp\left(C\left((\nu)^{1-1/n}+\ln(\nu)\right)\right),
\end{equation}
where $\nu\in{\Bbb N}.$

Let $\theta_1\in\{1,\,\ldots\,s\}. $Replacing $m$  by some bigger $m$
 and applying to $\beta^\wedge_{\theta_1}$ in (\ref{eq:xz2cdj})
 the Theorem 7, we see that
 there exist $A^\vee>0,\,m\in{\Bbb N},\,\alpha^\vee(\nu)>0$ with
 $\nu \in {m+{\Bbb N}-1}$  and the subspaces
 $R_{m,\theta_1,\theta_1}^\vee,\,\ldots,\,R_{m,\theta_1,s+1}^\vee$
 of the space $R_m=\Ker(\beta^\wedge_{\theta_1})$
 such that
$$\lim\limits_{\nu\to\infty}\alpha^\vee(\nu)=0,$$
$$R_m=R_{m,\theta_1,\theta_1}^\vee\oplus\,\ldots,\oplus\,
R_{m,\theta_1,s+1}^\vee,\,
\dim_{\Bbb C}(R_{m,\theta_1,i}^\vee)=e_i,\;\theta_1,\leq i\leq s+1,$$
if
\begin{equation}\label{eq:zcdd5a}
r\in R_{m,\theta_1,\theta}^\vee
\end{equation}
for some $\theta\in\{\theta_1,,...,s\}$, then
\begin{equation}\label{eq:zdc1a}
\exp(-A^\vee(\ln(\nu)+\nu^{1-1/k_\theta}))
\omega_{q,r}(m)(\rho_\theta)^{\nu}\le
\omega_{q,r)}(\nu),
\end{equation}
where $\nu\in {\Bbb N}+m-1;$ moreover, the spaces
$$R_{m,\theta_1,j}^\wedge =
 R_{m,\theta_1,j}^\vee \oplus\,\ldots\,\oplus R_{m,\theta_1,s+1}^\vee,$$
where $j=\theta_1,,\,\ldots,\, s+1$
 and natural projections $\pi_{\theta_1,j}^\vee$
 of the space  $R_{m,\theta_1,j}^\wedge$
 onto the space $R_{m,\theta_1,j}^\vee,$ where $j=\theta_1,,\,\ldots,\, s,$
 have the folloving properties:
\vskip2pt
if
\begin{equation}\label{eq:zcdd5b}
r\in R_{m,\theta_1,\theta}^\wedge
\end{equation}
for some $\theta\in\{\theta_1,,...,s\}$, then
\begin{equation}\label{eq:zdc1b}
\omega_{q,r}(\nu)\le(\rho_\theta)^{\nu}
\exp(A^\vee(\ln(\nu)+\nu^{1-1/k_\theta}))\omega_{q,r}(m),
\end{equation}
\begin{equation}\label{eq:zdc1c}
\omega_{q,r}(\nu)\ge(\omega_{q,\pi_{\theta_1,\theta}^\vee(r)}(m)
-\alpha^\vee(\nu)\omega_{q,r)}(m))\times
\end{equation}
$$
(\rho_\theta)^{\nu}
\exp(-A^\vee(\ln(\nu)+\nu^{1-1/k_\theta})),
$$
where $\nu\in {\Bbb N}+m-1;$ if $e_{s+1}>0,$ and
\begin{equation}\label{eq:zcdd6}
r\in R_{m,\theta_1,s+1}^\vee\,(=R_{m,\theta_1,s+1}^\wedge),
\end{equation}
 then
\begin{equation}\label{eq:zdc2}
\vert r(\nu)\vert\le(A^\vee/\nu)^{\nu/k^\ast}\omega_{q,r}(m),\end{equation}
where $\nu\in {\Bbb N}+m-1.$
\vskip4pt

Let
\begin{equation}\label{eq:z2cdd}
\beta=\beta^\ast_1=\beta_1,\,\eta=\beta^\wedge_{2}.
\end{equation}
Then, in view of (\ref{eq:yz2cdj}), (\ref{eq:2cdj}) holds.

\bf Lemma 9. \it Let all the conditions of the
Theorem 7 are fulfilled.
For sufficient big $m\in\Bbb N$ there exists the splitting monomorphism
 $$\psi^\wedge_2\colon \Ker(\eta)\to \Ker(\eta\beta)$$
with following properties:

(a)
$$\Ker(\beta)=V_{m,1}^\vee,\psi^\wedge_2(\Ker(\eta))=V_{m,2}^\wedge,$$
where $V_{m,1}^\vee$ and $V_{m,2}^\wedge$ are defined in the assertion
 of the Theorem 7;

(b) the monomorphism $\psi^\wedge_2$ maps isomorphically the space
$R_{m,2,j}^\vee,$ where $j=2,,\,\ldots,\, s+1,$ onto
the space $V_{m,j}^\vee;$

(c) the map $\beta_1\psi^wedge_2$ coincides with
 identiti map $\Ker(\eta)\to \Ker(\eta);$

(d) the restriction
 of the map $\psi^wedge_2\beta_1$ on the space $V_{m,2}^\wedge$
 coincides with the identity map $V_{m,2}^\wedge\to V_{m,2}^\wedge;$

(e) let $I$ be the identiti map $V_m\to V_m;$ then
 the natural projection $\pi_1$ of the space $V_{m,1}^\wedge=V_m$
 onto $V_{m,1}^\vee$ coincides with restriction
 of the map $I-\psi^\wedge_2\beta_1$ on the space $V_{m,1}^\vee;$

(f) if $j\in[2,s]\cap{\Bbb N},$ then natural projection $\pi_j$ of
 the space $V_{m,j}^\wedge$ onto the space$V_{m,j}^\vee$ coincides with
 the restriction of the map $\psi^\wedge_2\pi_{2,j}^\vee\beta_1$
 on the space $V_{m,j}^\wedge$ and the projection $\pi_{2,j}^\vee$
 coincides with the restriction of the map $\beta_1\pi_j,\psi^\wedge_2$
 on the space $R_{2,j}^\wedge.$

\bf Proof. \rm See the proof of the Theorem 7 in [\ref{r:ci0}],
 especially the section 2 and section 3. $\blacksquare.$

\bf Remark 2. \rm I use this opportunity to make a corrections
 in [\ref{r:ci0}]. On the page 13, third line from the bottom
 must stand $\psi_m\pi^\vee_j\phi_m$ instead $\psi_m\pi^\vee_j;$
 on the second line from the bottom
 must stand $\pi^\vee_j=\phi_m\pi_j\psi_m$ instead $\pi^\vee_j=\phi_m\pi_j.$

\bf Lemma 10. \it Let all the conditions of the Lemma 9 are fulfilled.
 Then (for sufficient big $m$)
\begin{equation}\label{eq:2ba}
V^\ast_{m,\theta}=\Ker(\beta_{1,\theta}^\ast),
\end{equation}
where $\theta=1,\,\ldots,\,s.$

\bf Proof. \rm We use induction on $s$ For $s=1$ or $\theta=1$
 the assertion of the Lemma directly follows from the Lemma 9.

Let $s>1$ and assertion of the lemma is true for $s-1.$ Let $\theta>1.$
 Since the monomorphism $\psi^\wedge_2$ maps isomorphically
 the space $R_{m,2,j}^\vee,$ where $j=2,,\,\ldots,\,\theta,$ onto
 the space $V_{m,j}^\vee$ and  the map $\beta_1\psi^\wedge_2$ coincides with
 identity map $\Ker(\eta)\to \Ker(\eta),$ it follows that
 the restriction of the map $\beta_1$ on
 the space $V_{m,2}^\vee\oplus\ldots\oplus V_{m,\theta}^\vee$
 is an isomorphism of $V_{m,2}^\vee\oplus\ldots\oplus V_{m,\theta}^\vee$
 onto  $R_{m,2,\theta}^\ast.$ According to the inductive hypothesis,
$$R_{m,2,\theta}^\ast=Ker(\beta^\ast_{2,\theta}).$$
 Therefore, if $y\in V_{m,2}^\vee\oplus\ldots\oplus V_{m,\theta}^\vee$
 then $\beta_{1,\theta}^\ast(y)=\beta_{2,\theta}^\ast(\beta_1(y))=0;$
  moreover, since also $V_{m,1}^\vee= \Ker(\beta_1),$ it follows that
 $V_{m,\theta}^\ast\subset\Ker(\beta_{1,\theta}^\ast).$
 Since, according to the Theorem 7,
\begin{equation}\label{eq:2bb}
\dim_{\mathbb C}(V_{m,\theta})=\sum\limits_{i=1}^\theta e_i=
\end{equation}
$$\dim_{\mathbb C}(\Ker(\beta_{1,\theta}^\ast))=
n_\theta:=\deg(\beta_{1,\theta}^\ast),$$
it follows that (\ref{eq:2ba}) holds. $\blacksquare.$

\bf Proof of the Theorem 8. \rm In view of (\ref{eq:2ba}), (\ref{eq:2bb}),
 we apply to $\beta_{1,\theta}^\ast$ the Lemma 3.  Since
 $\theta$ plays the role of $s$ in the lemma 3 now and
 $$\max(k_1,\,\ldots,\,k_\theta)\le k=\max(k_1,\,\ldots,\,k_s),$$
it follows that there exist $A=A_\theta>0,\,m\in{\Bbb N}$ such that
\begin{equation}\label{eq:2bc}
\exp(-A(ln(\nu)+\nu^{1-1/k}))(\rho_s)^{\nu}\omega_{n_\theta,y}(m)
\le\end{equation}
$$\omega_{n_\theta,y}(\nu)\le\omega_{n,y}(\nu)$$
for any $y\in V^\ast_{m,\theta}.$
 Since both the functions $y\to\omega_{n_theta,y}(m),y\in V^\ast_{m,\theta},$
 and $y\to\omega_{n,y}(m),y\in V^\ast_{m,\theta}$ are two norms
 on the $n_\theta$-dimensional space $V^\ast_{m,\theta}$ there
 exists  $B=B_\theta>0$ such that
$$\omega_{n,y}(m)\exp(-B)\le\omega_{n_\theta,y}(m)$$
for any $y\in V^\ast_{m,\theta}.$ Then
\begin{equation}\label{eq:2bd}
\exp(-A(ln(\nu)+\nu^{1-1/k})-B)(\rho_\theta)^{\nu}\omega_{n,y}(m)
\le
\omega_{n,y}(\nu)\end{equation}
for any $y\in V^\ast_{m,\theta}.\,\blacksquare$
{\begin{center}\large\bf \S 3.
 Proof of the theorem 9.\end{center}}
\bf Theorem 9. \it Let for some $\theta\in\{1\,\ldots,\, s\}$ is
 given a linear map $\xi_\theta$ of the space $V_{m,\theta}^\ast$
 into $V_{m,\theta+1}^\wedge.$ Let $I^\ast_\theta$ is
 the identity map $V_{m,\theta}^\ast\to V_{m,\theta}^\ast$
 Then for
$$V=(I_\theta^\ast+\xi_\theta)(V_{m,\theta}^\ast)$$
holds the assertion of the Theorem 10.

\bf  Proof of the theorem 9. \rm Since $V_{m,\theta}^\ast$ is
 finite-dimensional linear space over $\Bbb C,$ it follows that
 the map $\xi_\theta$ is continuous. Therevore there exists $C>0$ such that
\begin{equation}\label{eq:3be}
\omega_{n_\theta,\xi_\theta(y)}(m)\le
\exp(C)\omega_{n_\theta,\xi_\theta(y)}(m)
\end{equation}
for any $y\in V^\ast_{m,\theta}.$ According to the Theorem 7 and Theorem 8,
 there exist a numbers $A>0,\,m\in{\Bbb N},$ such that,
 if $y\in V^\ast_{m,\theta},\,\nu\in m-1+\Bbb N,$ then
\begin{equation}\label{eq:3bf}
\omega_{n,\xi_\theta(y)}(\nu)\le
\exp(A(\ln(\nu)+\nu^{1-1/k_{\theta+1}}))(\rho_\theta)^\nu
\omega_{n,\xi_\theta(y)}(m)
\end{equation}
for $\theta<s,$
\begin{equation}\label{eq:3bg}
\omega_{n_\theta,\xi_\theta(y)}(m)=0
\end{equation}
for $\theta=s$ and $k^\ast=k_{s+1}=0,$
\begin{equation}\label{eq:3bh}
\vert \xi_\theta(y)(\nu)\vert\le(A/\nu)^{\nu/k^\ast}
\omega_{n,\xi_\theta(y)}(m)
\end{equation}
for $\theta=s$ and $k^\ast=k_{s+1}>0,$
\begin{equation}\label{eq:3bi}
\exp(-A(ln(\nu)+\nu^{1-1/k}))(\rho_\theta)^{\nu}\omega_{n,y}(m)
\le
\omega_{n,y}(\nu).
\end{equation}
Therefore
\begin{equation}\label{eq:3cj}
\omega_{n,y+\xi_\theta(y)}(\nu)\ge\omega_{n,y}(\nu)-
\omega_{n,\xi_\theta(y)}(\nu)\ge
\exp(-A(ln(\nu)+\nu^{1-1/k}))(\rho_\theta)^{\nu}\omega_{n,y}(m)
\end{equation}
for $\theta=s$ and $k^\ast=k_{s+1}=0,$
\begin{equation}\label{eq:3ca}
\omega_{n,y+\xi_\theta(y)}(\nu)\ge\omega_{n,y}(\nu)-
\omega_{n,\xi_\theta(y)}(\nu)\ge
\end{equation}
$$\exp(-A(ln(\nu)+\nu^{1-1/k}))(\rho_\theta)^{\nu}\omega_{n,y}(m)-$$
$$\vert \xi_\theta(y)(\nu)\vert\le(A/\nu)^{\nu/k^\ast}
\exp(C)\omega_{n,y}(m)=$$
$$\exp(-A(ln(\nu)+\nu^{1-1/k}))(\rho_\theta)^{\nu}\omega_{n,y}(m)\times$$
$$\left(1+\exp(A(ln(\nu)+\nu^{1-1/k}))(\rho_\theta)^{-\nu}
(A/\nu)^{\nu/k^\ast}\exp(C)\right)$$
for $\theta=s$ and $k^\ast=k_{s+1}>0,$
\begin{equation}\label{eq:3cb}
\omega_{n,y+\xi_\theta(y)}(\nu)\ge\omega_{n,y}(\nu)-
\omega_{n,\xi_\theta(y)}(\nu)\ge
\end{equation}
$$\exp(-A(ln(\nu)+\nu^{1-1/k}))(\rho_\theta)^{\nu}\omega_{n,y}(m)-$$
$$(A/\nu)^{\nu/k^\ast}
\exp(C)\omega_{n,y}(m)=$$
$$\exp(-A(ln(\nu)+\nu^{1-1/k}))(\rho_\theta)^{\nu}\omega_{n,y}(m)\times$$
$$\left(1-\exp(A(ln(\nu)+\nu^{1-1/k}))(\rho_\theta)^{-\nu}
(A/\nu)^{\nu/k^\ast}\exp(C)\right)$$
for $\theta=s$ and $k^\ast=k_{s+1}>0,$
\begin{equation}\label{eq:3cc}
\omega_{n,y+\xi_\theta(y)}(\nu)\ge\omega_{n,y}(\nu)-
\omega_{n,\xi_\theta(y)}(\nu)\ge
\end{equation}
$$\exp(-A(ln(\nu)+\nu^{1-1/k}))(\rho_\theta)^{\nu}\omega_{n,y}(m)-$$
$$
\exp(A(\ln(\nu)+\nu^{1-1/k}))(\rho_\theta)^{\nu}\exp(C)\omega_{n,y}(m)=
$$
$$
\exp(-A(ln(\nu)+\nu^{1-1/k}))(\rho_\theta)^{\nu}\omega_{n,y}(m)\times
$$
$$\left(1-\exp(2A(ln(\nu)+\nu^{1-1/k}))(\rho_{\theta+1}/\rho_\theta)^{\nu}
\exp(C)\right)$$
for $\theta<s. \blacksquare$

\bf Remark 3. \rm The value of $\omega_{n,y}(m)$ in (\ref{eq:za4})
 may be much bigger than the value of $\omega_{n\pi_{\theta(y)}}(m).$
 Therefore we cannot (as show a simple examples)
 for fixed value of $A$ in the inequality (\ref{eq:za5}) to replace
 the linear space $V$ of the Theorem 10
 by the set $V_m\diagdown V^\wedge_{m,\theta+1}.$
 But in (\ref{eq:3cb}) -- (\ref{eq:3cc}) $\omega_{n,y}(m)$
 is carried out the brackets, and the value in the brackets
 tends to $1,$ when $\nu\to+\infty,$ being greather than $1/2$ for
 $\nu\in m-1+\Bbb N,$ with sufficient big $m\in\Bbb N,$
 which depends only from the equation.
{\begin{center}\large\bf \S 4.
 Proof of the theorem 10.\end{center}}
\bf  Proof of the theorem 10. \rm Let $\pi^\ast$ and $\pi^{\wedge}$
 are restrictions on $V$ of the natural endomorphisms of
 the space $V_m$ onto respectively $V_{m,\theta}^\ast$
 and $V_{m,\theta+1}^\wedge$ and let $I_0$ be the identity map $V\to V.$ Then
\begin{equation}\label{eq:4cb}
 I_0=\pi^\ast+\pi^{\wedge}
\end{equation}
Since $\Ker(\pi^\ast)\subset V\cap V_{m,\theta+1}^\wedge=\{0\},$
 it follows that $\pi^\ast$ is an isomorphism of the space $V$
 onto linear subspace $V^\prime$ of the space $V_{m,\theta}^\ast.$
 Let $\tau$ be the inverse isomorphism to the isomorphism $\pi^\ast.$
 Then $\pi^\ast\tau$ is identity map $V^\prime\to V^\prime$ and
 it follows from (\ref{eq:4cb}) that
\begin{equation}\label{eq:4cc}
\tau=\pi^\ast\tau+\pi^\wedge\tau.
\end{equation}
Clearly, the linear map $\pi^\wedge\tau:\, V^\prime\to V_{m,\theta+1}^\wedge$
 have an extension
\begin{equation}\label{eq:4cd}
\xi_\theta:\, V_{m,\theta}^\ast\to V_{m,\theta+1}^\wedge.
\end{equation}
It follows from (\ref{eq:4cc}) that
$$V=\tau(V^\prime)\subset (I+\xi_\theta)V_{m,\theta}^\ast,$$
where $I$ is identity map $V_{m,\theta}\to V_{m,\theta}$ and
 $\xi_\theta$ is a linear map in (\ref{eq:4cd}).

So, the Theorem 10 is Corollary of the Theorem 9. $\blacksquare$
{\begin{center}\large\bf \S 5.
 The case of the general differrence equation of the Poincar\'e type.
\end{center}}
 If (\ref{eq:xzc}) and (\ref{eq:yzc}) hold instead of (\ref{eq:zc}),
 then making use of the above arguments we obtain the following
 chainges in the Theorem 7 and theorem 10.
 For any $\varepsilon>0$ there exists
 a constant$A^\wedge=A^\wedge(\varepsilon)>0$ such that,
 if $y\in V_{m,\theta}^\wedge$ with some $\theta\in\{1,...,s\}$, then
 (instead of (\ref{eq:za2}) the following inequality holds)
\begin{equation}\label{eq:yza2}
\omega_{n,y}(\nu)\le
\exp(A^\wedge(\rho_\theta\exp(\varepsilon))^{\nu}\omega_{n,y}(m),
\end{equation}
where $\nu\in m+{\Bbb N}-1;$
 if $y\in V_{m,s+1}^\vee\,(=V_{m,s+1}^\wedge)$, then
 (instead of (\ref{eq:xza3}) the following inequality holds)
\begin{equation}\label{eq:yza3}
\omega_{n,y}(\nu)\le
\exp(A^\wedge(\exp(-\varepsilon))^{\nu}\omega_{n,y}(m),
\end{equation}
where $\nu\in m+{\Bbb N}-1.$

Let further $V$ be an arbitrary linear subspace of $V_m$ such that
$$V\cap V_{m,\theta+1}=\{0\},$$
where $\theta\in\{1\,\ldots,\, s\}.$
 Then for this subspace $V$ and any $\varepsilon>0$
 there exists a constant $A^\vee=A^\vee(V,\varepsilon)>0$
 such that (instead of (\ref{eq:za5}) the following inequality holds)
\begin{equation}\label{eq:xza5}
\exp(-A^\vee)(\rho_\theta\exp(-\varepsilon))^{\nu}\omega_n(y)(m)
\le
\omega_{n,y}(\nu)\end{equation}
where $y\in V$ and $\nu\in m+{\Bbb N}-1.$

\bf Corollary. \rm (See [\ref{r:bj}], Theorem 3 and [\ref{r:ch0}], Lemma 16).
 \it Let as before (\ref{eq:xzc}) and (\ref{eq:yzc}) hold
 instead of (\ref{eq:zc}). Let $V$ be a $r$-dimensional subspace of $V_m$,
 let
$$V\cap V^\vee_{m,s+1}=\{0\}$$
 and let
 $\{y_1(\nu)\,,\ldots\,,y_r(\nu)\}$ be  a basis  of the space $V$. Let
$$k_3(V)=\max\{k\in{\Bbb Z}\colon 1\leq k\leq s,\;
V\subset V^\wedge_{m,k}\},$$
and
$$k_4(V)=\min\{k\in {\Bbb Z}\colon 1\leq k\leq s,\;
V\cap V^\wedge_{m,k+1}=\{0\}\}.$$
For $X = (x_1\,,\ldots\,,x_r),\;X\in{\Bbb C}^r$, let
 $$h(X)=max\{\vert x_1\vert\,,\ldots\,,\vert x_r\vert\},$$
 $$y=y^\vee(X,\nu)=x_1 y^\vee_1(\nu)\,+\ldots\,+ x_ry^\vee_r(\nu).$$
 Then for every $\varepsilon\in(0,1)$ there exist
 $C_3(\varepsilon)>0$ and $C_4(\varepsilon)>0$ such that
$$
C_4(\varepsilon)(\rho_{k_4}(1-\varepsilon))^\nu h(X)\le\omega_{n,y}(\nu)\le
C_3(\varepsilon) (\rho_{k_3 }+\varepsilon)^\nu h(X).$$

\bf Proof \rm The functions $h(X)$ and the restriction
 of $\omega_{n,y}(m)$ on $V$ are two norms on the $r$-dimensional
 over $\Bbb C$ linear space $V.$ Therefore there exists a constant $C_5>0$
 such that $ h(X)\le C_5\omega_{n,y}(m)$ and $\omega_{n,y}(m)\le C_5h(X)$
 The assertion of the Corollary directly follows
 from (\ref{eq:yza2}) and (\ref{eq:xza5}) now. $\blacksquare$
{\begin{center}\large\bf References.\end{center}}
\footnotesize
\vskip4pt
\refstepcounter{r}\noindent[\ther]
R.Ap\'ery, Interpolation des fractions continues\\
\hspace*{3cm}  et irrationalite de certaines constantes,\\
\hspace*{3cm} Bulletin de la section des sciences du C.T.H., 1981, No 3,
 37 -- 53;
\label{r:cd}\\
\refstepcounter{r}
\noindent[\ther]
F.Beukers, A note on the irrationality of $\zeta (2)$ and $\zeta (3),$\\
\hspace*{3cm} Bull. London Math. Soc., 1979, 11, 268 --  272;
\label{r:ce}\\
\refstepcounter{r}
\noindent[\ther]
 A.van der Porten,
 A proof that Euler missed...Ap\'ery's proof of the irrationality of
 $\zeta (3),$\\
\hspace*{3cm} Math Intellegencer, 1979, 1, 195 -- 203;\label{r:cf}\\
\refstepcounter{r}
\noindent[\ther]
 W. Maier, Potenzreihen irrationalen Grenzwertes,\\
\hspace*{3cm} J.reine angew. Math.,
 156, 1927, 93 -- 148;\label{r:cg}\\
\refstepcounter{r}\noindent[\ther]
 E.M. Niki\u sin, On irrationality of the values of the functions F(x,s)
 (in Russian),\\
\hspace*{3cm} Mat.Sb. 109 (1979), 410 -- 417;\\
\hspace*{3cm}
 English transl. in Math. USSR Sb. 37 (1980), 381 -- 388;\label{r:ch}\\
\refstepcounter{r}\noindent[\ther]
 G.V. Chudnovsky, Pade approximations to the generalized
 hyper-geometric functions\\
\hspace*{3.4cm}  I,J.Math.Pures Appl., 58, 1979,
 445 -- 476;\label{r:dj}\\
\refstepcounter{r}\noindent[\ther]
\rule{2.7cm}{.3pt},
 Transcendental numbers,  Number Theory,Carbondale,\\
\hspace*{3.4cm}  Lecture Notes in Math, Springer-Verlag, 1979, 751, 45 -- 69;
\label{r:da}\\
\refstepcounter{r}\noindent[\ther]
\rule{2.7cm}{.3pt}, Approximations rationelles des logarithmes
 de nombres rationelles\\
\hspace*{3.4cm} C.R.Acad.Sc. Paris, S\'erie A, 1979, 288, 607 -- 609;
\label{r:db}\\
\refstepcounter{r}\noindent[\ther]
\rule{2.7cm}{.3pt}, Formules d'Hermite pour les approximants de Pad\'e de
logarithmes\\
\hspace*{3.4cm} et de fonctions bin\^omes, et mesures d'irrationalit\'e,\\
\hspace*{3.4cm} C.R.Acad.Sc. Paris, S\'erie A, 1979, t.288, 965 -- 967;
\label{r:dc}\\
\refstepcounter{r}\noindent[\ther]
\rule{2.5cm}{.3pt},Un syst\'me explicite d'approximants de Pad\'e\\
\hspace*{3.4cm} pour les fonctions hyp\'erg\'eometriques g\'en\'eralies\'ees,\\
\hspace*{3.4cm} avec applications a l'arithm\'etique,\\
\hspace*{3.4cm}
 C.R.Acad.Sc. Paris, S\'erie A, 1979, t.288, 1001 -- 1004;\label{r:dd}\\
\refstepcounter{r}\noindent[\ther]
\rule{2.5cm}{.3pt},
 Recurrenses defining Rational Approximations\\
\hspace*{3.4cm} to the irrational numbers, Proceedings\\
\hspace*{3.4cm}  of the Japan Academie,
 Ser. A, 1982, 58, 129 -- 133;
\label{r:de}\\
\refstepcounter{r}\noindent[\ther]
\rule{2.5cm}{.3pt}, On the method of Thue-Siegel,\\
\hspace*{3.4cm} Annals of Mathematics, 117 (1983), 325 -- 382;
\label{r:df}\\
\refstepcounter{r}\noindent[\ther]
K.Alladi  and M. Robinson, Legendre polinomials and irrationality,\\
\hspace*{6cm}J. Reine Angew.Math., 1980, 318, 137 -- 155;
\label{r:dg}\\
\refstepcounter{r}\noindent[\ther]
A. Dubitskas, An approximation of logarithms of some numbers,\\
\hspace*{3cm} Diophantine approximations II,Moscow, 1986, 20 -- 34;
\label{r:dh}\\
\refstepcounter{r}\noindent [\ther] \rule{2cm}{.3pt},
 On approximation of $\pi/ \sqrt {3}$ by rational fractions,\\
\hspace*{3cm} Vestnik MGU, series 1, 1987, 6, 73 -- 76;
\label{r:ej}\\
\refstepcounter{r}\noindent[\ther]
 S.Eckmann, \"Uber die lineare Unaqbhangigkeit der Werte gewisser Reihen,\\
\hspace*{3cm} Results in Mathematics, 11, 1987, 7 -- 43;
\label{r:ea}\\
\refstepcounter{r}\noindent[\ther]
 M.Hata, Legendre type polinomials and irrationality mesures,\\
\hspace*{3cm} J. Reine Angew. Math., 1990, 407, 99 -- 125;
\label{r:eb}\\
\refstepcounter{r}\noindent[\ther]
 A.O. Gelfond, Transcendental and algebraic numbers (in Russian),\\
\hspace*{3cm} GIT-TL, Moscow, 1952;
\label{r:ec}\\
\refstepcounter{r}\noindent[\ther]
H.Bateman and A.Erd\'elyi, Higher transcendental functions,1953,\\
\hspace*{3cm}  New-York -- Toronto -- London,
 Mc. Grow-Hill Book Company, Inc.;
\label{r:ed}\\
\refstepcounter{r}\noindent[\ther]
 E.C.Titchmarsh, The Theory of Functions, 1939, Oxford University Press;
\label{r:ee}\\
\refstepcounter{r}\noindent[\ther]
E.T.Whittaker and G.N. Watson, A course of modern analysis,\\
\hspace*{3cm} 1927, Cambridge University Press;
\label{r:ef}\\
\refstepcounter{r}\noindent[\ther]
 O.Perron, \"Uber die Poincaresche Differenzengleichumg,\\
\hspace*{3cm} Journal f\"ur die reine und angewandte mathematik,\\
\hspace*{3cm} 1910, 137,  6 -- 64;\label{r:a}\\
\refstepcounter{r}
\noindent [\ther] A.O.Gelfond, Differenzenrechnung (in Russian),
 1967, Nauka, Moscow.
\label{r:b}\\
\refstepcounter{r}\noindent [\ther]
 A.O.Gelfond and I.M.Kubenskaya, On the theorem of Perron\\
\hspace*{6cm} in the theory of differrence equations  (in Russian),\\
\hspace*{3cm} IAN USSR, math. ser., 1953, 17, 2,  83 --  86.
\label{r:c}\\
\refstepcounter{r}
\noindent [\ther] M.A.Evgrafov, New proof of the theorem of Perron\\
\hspace*{3cm} (in Russian),IAN USSR, math. ser., 1953, 17, 2, 77 --  82;
\label{r:d}\\
\refstepcounter{r}
\noindent [\ther] G.A.Frejman, On theorems of of Poincar\'e and Perron\\
\hspace*{3cm} (in Russian), UMN, 1957, 12, 3 (75), 243 -- 245;
\label{r:e}\\
\refstepcounter{r}
\noindent [\ther]  N.E.N\"orlund, Differenzenrechnung, Berlin,
 Springer Verlag, 1924;\label{r:f}\\
\refstepcounter{r}
\noindent [\ther] I.M.Vinogradov, Foundtions of the Number Theory,
(in Russian), 1952, GIT-TL;\label{r:g}\\
\refstepcounter{r}
\noindent [\ther] \rule{2cm}{.3pt},
J. Diedonne, Foundations of modern analysis,\\
\hspace*{3cm}
 Institut des Hautes \'Etudes Scientifiques, Paris,\\
\hspace*{3cm} Academic Press, New York and London, 1960\label{r:ei}\\
\refstepcounter{r}
\noindent [\ther] CH.-J. de la Vall\'ee Poussin, Course d'analyse
infinit\'esimale,\\
\hspace*{3cm} Russian translation by G.M.Fikhtengolts,\\
\hspace*{3cm} GT-TI, 1933;
\label{r:cb0}\\
\refstepcounter{r}
\noindent [\ther] H.WEYL, Algebraic theory of numbers,\\
\hspace*{3cm}  1940,\\
\hspace*{3cm} Russian translation by L.I.Kopejkina;
\label{r:cb1}\\
\refstepcounter{r}
\noindent [\ther] G.Rhin and C. Viola.\\
\hspace*{3cm} On a prmutation group related to $\zeta(2),$\\
\hspace*{3cm} Acta Arithmetica 77(1996), 25 -- 56;
\label{r:rv1}\\
\refstepcounter{r}
\noindent [\ther] G.Rhin and C. Viola. The group structure for $\zeta(3),$\\
\hspace*{3cm} Acta Arithmetica 97(2001), 269 -- 293;
\label{r:rv2}\\
\refstepcounter{r}
\noindent [\ther] M. Hata. A new irrationality measure for$\zeta(3),$\\
\hspace*{3cm} Acta Arithmetica 92(2000), 47 -- 57;
\label{r:h2}\\
\refstepcounter{r}
\noindent [\ther] T.Rivual. La focntion z\^eta de Riemann prend une infinit\'e\\
\hspace*{3cm} valeurs irrationnelles aux entires impairs,\\
\hspace*{3cm} C.R.Acad.Sci.Paris, s\'erie 1, p. 267 -- 270, 2000.
\label{r:riv1}\\
\refstepcounter{r}
\noindent [\ther] K.Boll, T.Rivual. Irrationalit\'e d`une infinit\'e\\
\hspace*{3cm} valeurs la focntion z\^eta aux entires impairs,\\
\hspace*{3cm} Invent. math., 146, 193 -- 207 (2001);
\label{r:bolriv1}\\
\refstepcounter{r}
\noindent [\ther] W.Zudilin. One from the nubers
$\zeta(5),\, \zeta(7),\, \zeta(9),\, \zeta(11)$ is irrational.\\
\hspace*{3cm} (in Russian) UMN 56(4) (2001), 149 -- 150;
\label{r:z1}\\
\refstepcounter{r}
\noindent [\ther] W.Zudilin. Ap\'ery's Theorem and problems for the values\\
\hspace*{3cm} of $\zeta$-function of Riemann and their $q-$analogs.\\
\hspace*{3cm} (in Russian) Mocow state university, doctoral thesis,\\
\hspace*{3cm} Mocow 2004;
\label{r:z2}\\
\noindent [\ther] L.A.Gutnik, On the decomposition
 of the difference operators of Poincar\'e type\\
\hspace*{3cm}  (in Russian),
VINITI, Moscow, 1992, 2468 -- 92, 1 -- 55;
\label{r:h}\\
\refstepcounter{r}
\noindent [\ther] \rule{2cm}{.3pt},  On the decomposition
 of the difference operators\\
\hspace*{3cm}  of Poincar\'e type in Banach algebras\\
\hspace*{3cm} (in Russian), VINITI, Moscow, 1992, 3443 -- 92, 1 -- 36;
\label{r:i}\\
\refstepcounter{r}
\noindent [\ther] \rule{2cm}{.3pt}, On the difference equations of Poincar\'e
 type\\
\hspace*{3cm}  (in Russian), VINITI, Moscow 1993, 443 -- B93, 1 -- 41;
\label{r:aj}\\
\refstepcounter{r}
\noindent [\ther] \rule{2cm}{.3pt}, On the difference equations
 of Poincar\'e type in normed algebras\\
\hspace*{3cm} (in Russian), VINITI, Moscow, 1994, 668 -- B94, 1 -- 44;
\label{r:aa}\\
\refstepcounter{r}
\noindent [\ther] \rule{2cm}{.3pt}, On the decomposition of
 the difference equations  of Poincar\'e type\\
\hspace*{3cm} (in Russian),
 VINITI, Moscow, 1997, 2062 -- B97, 1 -- 41;
\label{r:ab}\\
\refstepcounter{r}
\noindent [\ther] \rule{2cm}{.3pt}, The difference equations
 of Poincar\'e type\\
\hspace*{3cm} with characteristic polynomial having roots equal to zero\\
\hspace*{3cm}(in Russian), VINITI, Moscow, 1997, 2418 -- 97, 1 -- 20;
\label{r:ac}\\
\refstepcounter{r}
\noindent [\ther]  \rule{2cm}{.3pt}, On the behavior of solutions\\
\hspace*{3cm} of difference equations of Poincar\'e type\\
\hspace*{3cm}
 (in Russian), VINITI, Moscow, 1997, 3384 -- B97, 1 - 41;
\label{r:ad}\\
\refstepcounter{r}
\noindent [\ther] \rule{2cm}{.3pt}, On the variability of solutions
of difference equations of Poincar\'e type\\
\hspace*{3cm}
 (in Russian),  VINITI, Moscow, 1999, 361 -- B99, 1 -- 9;
\label{r:ae}\\
\refstepcounter{r}
\noindent [\ther] \rule{2cm}{.3pt}, To the question of the variability
 of solutions\\
\hspace*{3cm} of difference equations of Poincar\'e type (in Russian),\\
\hspace*{3cm} VINITI, Moscow, 2000, 2416 -- B00, 1 -- 22;\label{r:af}\\
\refstepcounter{r}
\noindent [\ther] \rule{2cm}{.3pt}, On linear forms with coefficients
 in ${\mathbb N}\zeta(1+\mathbb N),$\\
\hspace*{3cm} Max-Plank-Institut f\"ur Mathematik,\\
\hspace*{3cm} Bonn, Preprint Series, 2000, 3, 1 -- 13;
\label{r:ag}\\
\refstepcounter{r}
\noindent [\ther] \rule{2cm}{.3pt},
On the Irrationality of Some Quantyties Containing $\zeta (3)$ (in Russian),\\
\hspace*{3cm} Uspekhi Mat. Nauk, 1979, 34, 3(207), 190;
\label{r:di}\\
\refstepcounter{r}
\noindent [\ther] \rule{2cm}{.3pt}, On the Irrationality of Some Quantities
 Containing $\zeta (3),$\\
\hspace*{3cm} Eleven papers translated from the Russian,\\
\hspace*{3cm} American Mathematical Society, 1988, 140, 45 - 56;
\label{r:ah}\\
\refstepcounter{r}
\noindent [\ther] \rule{2cm}{.3pt}, Linear independence over $\mathbb Q$
 of dilogarithms at rational points\\
\hspace*{3cm} (in Russian), UMN, 37 (1982), 179-180;\\
\hspace*{3cm}english transl. in Russ. Math. surveys 37 (1982), 176-177;
\label{r:ai}\\
\refstepcounter{r}
\noindent [\ther] \rule{2cm}{.3pt}, On a measure of the irrationality
 of dilogarithms at rational points\\
\hspace*{3cm} (in Russian), VINITI, 1984, 4345-84, 1 -- 74;
\label{r:bj}\\
\refstepcounter{r}
\noindent [\ther] \rule{2cm}{.3pt}, To the question of the smallness
of some linear forms\\
\hspace*{3cm} (in Russian), VINITI, 1993, 2413-B93, 1 -- 94;
\label{r:ba}\\
\refstepcounter{r}
\noindent [\ther] \rule{2cm}{.3pt}, About linear forms, whose coefficients
 are logarithms\\
\hspace*{3cm} of algebraic numbers  (in Russian),\\
\hspace*{3cm} VINITI, 1995, 135-B95, 1 -- 149;
\label{r:bb}\\
\refstepcounter{r}
\noindent [\ther] \rule{2cm}{.3pt}, About systems of vectors, whose
 coordinates\\
\hspace*{3cm} are linear combinations of logarithms of algebraic numbers\\
\hspace*{3cm} with algebraic coefficients  (in Russian),\\
\hspace*{3cm} VINITI, 1994, 3122-B94, 1 -- 158;
\label{r:bc}\\
\refstepcounter{r}
\noindent [\ther] \rule{2cm}{.3pt}, On the linear forms, whose\\
\hspace*{3cm} coefficients are $\mathbb A$ - linear combinations\\
\hspace*{3cm}  of logarithms of $\mathbb A$ - numbers,\\
\hspace*{3cm} VINITI, 1996,  1617-B96, pp. 1 -- 23.
\label{r:bc1}\\
\refstepcounter{r}
\noindent [\ther] \rule{2cm}{.3pt}, On systems of linear forms, whose\\
\hspace*{3cm} coefficients are $\mathbb A$ - linear combinations\\
\hspace*{3cm}  of logarithms of $\mathbb A$ - numbers,\\
\hspace*{3cm} VINITI, 1996,  2663-B96, pp. 1 -- 18.
\label{r:bc2}\\
\refstepcounter{r}
\noindent [\ther] \rule{2cm}{.3pt}, About linear forms, whose coefficients\\
\hspace*{3cm} are $\mathbb Q$-proportional to the number $\log 2, $
 and the values\\
\hspace*{3cm} of $\zeta (s)$ for integer $s$ (in Russian),\\
\hspace*{3cm} VINITI, 1996, 3258-B96, 1 -- 70;
\label{r:bd}\\
\refstepcounter{r}
\noindent [\ther] \rule{2cm}{.3pt}, The lower estimate for some linear forms,\\
\hspace*{3cm} coefficients of which are proportional to the values\\
\hspace*{3cm} of $\zeta (s)$ for integer $s$ (in Russian),\\
\hspace*{3cm} VINITI, 1997, 3072-B97, 1 -- 77;
\label{r:be}\\
\refstepcounter{r}
\noindent [\ther] \rule{2cm}{.3pt},
 On linear forms with coefficients in
${\mathbb N} \zeta(1 + {\mathbb N}) $\\
\hspace*{3cm} Max-Plank-Institut f\"ur Mathematik, Bonn,\\
\hspace*{3cm} Preprint Series, 2000, 3, 1 -- 13;\label{r:cc0}\\
\refstepcounter{r}
\noindent [\ther] \rule{2cm}{.3pt},
 On linear forms with coefficients
 in ${\mathbb N}\zeta(1+\mathbb N)$\\
\hspace*{3cm} (the detailed version,part 1),
 Max-Plank-Institut f\"ur Mathematik,\\
\hspace*{3cm} Bonn, Preprint Series, 2001, 15, 1 -- 20;
\label{r:bf}\\
\refstepcounter{r}
\noindent [\ther] \rule{2cm}{.3pt}, On linear forms with coefficients
 in ${\mathbb N}\zeta(1+\mathbb N)$\\
\hspace*{3cm} (the detailed version,part 2),
 Max-Plank-Institut f\"ur Mathematik,\\
\hspace*{3cm} Bonn, Preprint Series,2001, 104, 1 -- 36;
\label{r:bg}\\
\refstepcounter{r}
\noindent [\ther] \rule{2cm}{.3pt}, On linear forms with coefficients
in ${\mathbb N}\zeta(1+{\mathbb N})$\\
\hspace*{3cm} (the detailed version,part 3),
 Max-Plank-Institut f\"ur Mathematik,\\
\hspace*{3cm} Bonn, Preprint Series, 2002, 57, 1 -- 33;
\label{r:bh}\\
\refstepcounter{r}
\noindent [\ther] \rule{2cm}{.3pt},
On the rank over ${\mathbb Q}$ of some real matrices (in Russian),\\
\hspace*{3cm} VINITI, 1984, 5736-84; 1 -- 29;\label{r:eg}\\
\refstepcounter{r}
\noindent [\ther] \rule{2cm}{.3pt},
On the rank over ${{\mathbb Q}}$ of some real matrices,\\
\hspace*{3cm} Max-Plank-Institut f\"ur Mathematik,\\
\hspace*{3cm} Bonn, Preprint Series, 2002, 27, 1 -- 32;\label{r:eh}\\
\refstepcounter{r}
\noindent [\ther] \rule{2cm}{.3pt}, On linear forms with coefficients
 in ${\mathbb N}\zeta(1+{\mathbb N})$\\
\hspace*{3cm} (the detailed version, part 4),
 Max-Plank-Institut f\"ur Mathematik,\\
\hspace*{3cm} Bonn, Preprint Series, 2002, 142, 1 -- 27;
\label{r:bi}\\
\refstepcounter{r}
\noindent [\ther] \rule{2cm}{.3pt}, On the dimension of some linear spaces\\
\hspace*{3cm} over finite extension of ${\mathbb Q}$ (part 2),\\
\hspace*{3cm} Max-Plank-Institut f\"ur Mathematik, Bonn, Preprint Series,\\
\hspace*{3cm} 2002, 107, 1 -- 37;
\label{r:cj}\\
\refstepcounter{r}
\noindent [\ther] \rule{2cm}{.3pt}, On the dimension of some linear spaces
 over $\mathbb Q$ (part 3),\\
\hspace*{3cm} Max-Plank-Institut f\"ur Mathematik, Bonn,\\
\hspace*{3cm} Preprint Series, 2003, 16, 1 -- 45.
\label{r:ca}\\
\refstepcounter{r}
\noindent [\ther] \rule{2cm}{.3pt},On the difference equation
 of Poincar\'e type (Part 1).\\
\hspace*{3cm} Max-Plank-Institut f\"ur Mathematik, Bonn,\\
\hspace*{3cm} Preprint Series, 2003, 52, 1 -- 44.
\label{r:cb}\\
\refstepcounter{r}
\noindent [\ther] \rule{2cm}{.3pt},  On the dimension of some linear
 spaces over $\mathbb Q,$ (part 4)\\
\hspace*{3cm} Max-Plank-Institut f\"ur Mathematik, Bonn,\\
\hspace*{3cm} Preprint Series, 2003, 73, 1 -- 38.
\label{r:cc}\\
\refstepcounter{r}
\noindent [\ther] \rule{2cm}{.3pt}, On linear forms with coefficients
 in ${\mathbb N}\zeta(1+\mathbb N)$\\
\hspace*{3cm} (the detailed version, part 5),\\
\hspace*{3cm} Max-Plank-Institut f\"ur Mathematik, Bonn,\\
\hspace*{3cm} Preprint Series, 2003, 83, 1 -- 13.
\label{r:cd0}\\
\refstepcounter{r}
\noindent [\ther] \rule{2cm}{.3pt}, On linear forms with coefficients
 in ${\mathbb N}\zeta(1+\mathbb N)$\\
\hspace*{3cm} (the detailed version, part 6),\\
\hspace*{3cm} Max-Plank-Institut f\"ur Mathematik, Bonn,\\
\hspace*{3cm} Preprint Series, 2003, 99, 1 -- 33.
\label{r:ce0}\\
\refstepcounter{r}
\noindent [\ther] \rule{2cm}{.3pt},On the difference equation
of Poincar\'e type (Part 2).\\
\hspace*{3cm} Max-Plank-Institut f\"ur Mathematik, Bonn,\\
\hspace*{3cm} Preprint Series, 2003, 107, 1 -- 25.
\label{r:cf0}\\
\refstepcounter{r}
\noindent [\ther] \rule{2cm}{.3pt},
On the asymptotic behavior of solutions\\
\hspace*{3cm}  of difference equation (in English).\\
\hspace*{3cm} Chebyshevskij sbornik,
 2003, v.4, issue 2, 142 -- 153.
\label{r:cg0}\\
\refstepcounter{r}
\noindent [\ther] \rule{2cm}{.3pt}, On linear combinations of logarithms\\
\hspace*{3cm} of algebraic  numbers with algebraic coefficients.\\
\hspace*{3cm}  International conference\\
\hspace*{3cm}  "Diophantine analysis, uniform distributions\\
\hspace*{3cm}  and applications"\\
\hspace*{3cm}  August 25-30, 2003, Minsk, Belarus\\
\hspace*{3cm} {\footnotesize Abstracts}\\
\hspace*{3cm}  pp. 16 -- 17.
\label{r:ch1}\\
\refstepcounter{r}
\noindent [\ther] \rule{2cm}{.3pt}, On linear forms with coefficients
 in ${\mathbb N}\zeta(1+\mathbb N),$\\
\hspace*{3cm} Bonner Mathematishe Schriften Nr. 360,\\
\hspace*{3cm} Bonn, 2003, 360.
\label{r:ch0}\\
\refstepcounter{r}
\noindent [\ther] \rule{2cm}{.3pt}, On linear forms with coefficients
 in ${\mathbb N}\zeta(1+\mathbb N)$\\
\hspace*{3cm} (the detailed version, part 7),\\
\hspace*{3cm} Max-Plank-Institut f\"ur Mathematik, Bonn,\\
\hspace*{3cm} Preprint Series, 2004, 88, 1 -- 27.
\label{r:zch0}\\
\refstepcounter{r}
\noindent [\ther] \rule{2cm}{.3pt},On the difference equation
of Poincar\'e type (Part 3).\\
\hspace*{3cm} Max-Plank-Institut f\"ur Mathematik, Bonn,\\
\hspace*{3cm} Preprint Series, 2004, 09, 1 -- 34.
\label{r:ci0}\\
\refstepcounter{r}
\noindent [\ther] \rule{2cm}{.3pt},  On the dimension of some linear
 spaces over $\mathbb Q,$ (part 5)\\
\hspace*{3cm} Max-Plank-Institut f\"ur Mathematik, Bonn,\\
\hspace*{3cm} Preprint Series, 2004,1 -- 42.
\label{r:ci1}
\vskip 10pt
 {\it E-mail:}{\sl\ gutnik$@@$gutnik.mccme.ru}
\end{document}